\documentclass[12pt]{amsart}

\usepackage{a4wide}

\usepackage{enumerate}

\usepackage{tikz}

\usetikzlibrary{arrows}

\newtheorem{theorem}{Theorem}[section]

\theoremstyle{definition}
\newtheorem{remark}[theorem]{Remark}
\newtheorem{example}[theorem]{Example}

\newcommand{\C}{\ensuremath{\mathbb{C}}}
\newcommand{\R}{\ensuremath{\mathbb{R}}}
\renewcommand{\H}{\ensuremath{\mathbb{H}}}
\renewcommand{\Im}{\ensuremath{\mathrm{Im}\;}}

\renewcommand{\P}{\ensuremath{\mathsf{P}}}
\newcommand{\Hy}{\ensuremath{\mathsf{H}}}
\newcommand{\SO}{\ensuremath{\mathsf{SO}}}
\newcommand{\SU}{\ensuremath{\mathsf{SU}}}
\newcommand{\SL}{\ensuremath{\mathsf{SL}}}
\newcommand{\Sg}{\ensuremath{\mathsf{S}}}
\newcommand{\U}{\ensuremath{\mathsf{U}}}
\newcommand{\Sp}{\ensuremath{\mathsf{Sp}}}
\newcommand{\Spin}{\ensuremath{\mathsf{Spin}}}
\newcommand{\g}[1]{\ensuremath{\mathfrak{#1}}}
\newcommand{\cal}[1]{\ensuremath{\mathcal{#1}}}

\newcommand{\Ss}{\ensuremath{\mathcal{S}}}
\newcommand{\Cl}{\ensuremath{\mathrm{Cl}}}

\newcommand{\End}{\ensuremath{\mathrm{End}}}
\newcommand{\rank}{\ensuremath{\mathrm{rank}\;}}

\DeclareMathOperator{\spann}{span}

\DeclareMathOperator{\tr}{tr}

\DeclareMathOperator{\Ad}{Ad}

\DeclareMathOperator{\Isom}{Isom}
\DeclareMathOperator{\ad}{ad}
\DeclareMathOperator{\id}{id}
\DeclareMathOperator{\Id}{Id}

\begin{document}

\title[Submanifold geometry in symmetric spaces of noncompact type]{Submanifold geometry in symmetric spaces\\of noncompact type}

\author[J.~C.~D\'iaz-Ramos]{J.~Carlos~D\'iaz-Ramos}
\address{Department of Mathematics, University of Santiago de Compostela, Spain.}
\email{josecarlos.diaz@usc.es}
\author[M.~Dom\'{\i}nguez-V\'{a}zquez]{Miguel Dom\'{\i}nguez-V\'{a}zquez}
\address{Instituto de Ciencias Matem\'aticas (CSIC-UAM-UC3M-UCM), Madrid, Spain.}
\email{miguel.dominguez@icmat.es}
\author[V.~Sanmart\'in-L\'opez]{V\'ictor Sanmart\'in-L\'opez}
\address{Department of Mathematics, University of Santiago de Compostela, Spain.}
\email{victor.sanmartin@usc.es}

\thanks{The authors have been supported by the project MTM2016-75897-P (AEI/FEDER, Spain). The second author acknowledges support by the ICMAT Severo Ochoa project SEV-2015-0554 (MINECO, Spain), as well as by the European Union's Horizon 2020 research and innovation programme under the Marie Sk\l{}odowska-Curie grant agreement No.~745722.}



\begin{abstract}
In this survey article we provide an introduction to submanifold geometry in symmetric spaces of noncompact type. We focus on the construction of examples and the classification problems of homogeneous and isoparametric hypersurfaces, polar and hyperpolar actions, and homogeneous CPC submanifolds.
\end{abstract}


\maketitle

\section{Introduction}
According to the original definition given by Cartan~\cite{Ca26}, a Riemannian symmetric space is a Riemannian manifold characterized by the property that curvature is invariant under parallel translation. This geometric definition has the surprising effect of bringing the theory of Lie groups into the picture, and it turns out that Riemannian symmetric spaces are intimately related to semisimple Lie groups. To a large extent, many geometric problems in symmetric spaces can be reduced to the study of properties of semisimple Lie algebras, thus transforming difficult geometric questions into linear algebra problems that one might be able to solve.

For this reason, the family of Riemannian symmetric spaces has been a setting where many geometric properties can be tackled and tested. They are often a source of examples and counterexamples. The set of symmetric spaces is a large family encompassing many interesting examples of Riemannian manifolds such as spaces of constant curvature, projective and hyperbolic spaces, Grassmannians, compact Lie groups and more. Apart from Differential Geometry, symmetric spaces have also been studied from the point of view of Global Analysis and Harmonic Analysis, being noncompact symmetric spaces of particular relevance (see for example~\cite{He08}).  They are also an outstanding family in the theory of holonomy, constituting a class of their own in Berger's classification of holonomy groups.

Our interest in symmetric spaces comes from a very general question: the relation between symmetry and shape. In a broad sense, the symmetries of a mathematical object are the transformations of that object that leave it invariant. These symmetries impose several constraints that reduce the degrees of freedom of the object, and imply a regularity on its shape.  More concretely, in Submanifold Geometry of Riemannian manifolds, our symmetric objects will actually be (extrinsically) {homogeneous submanifolds}, that is, submanifolds of a given Riemannian manifold that are orbits of a subgroup of isometries of the ambient manifold.  In other words, a submanifold $P$ of a Riemannian manifold $M$ is said to be homogeneous if for any two points $p$, $q\in P$ there exists and isometry $\varphi$ of $M$ such that $\varphi(P)=P$ and $\varphi(p)=q$. The symmetries of $M$ are precisely the isometries $\varphi$ in this definition.  Therefore, the study of homogeneous submanifolds makes sense only in ambient manifolds with a large group of isometries, and thus, the class of Riemannian symmetric spaces is an ideal setup for this problem.

Roughly speaking (see Section~\ref{sec:symmetric}) there are three types of symmetric spaces: Euclidean spaces, symmetric spaces of compact type (in case the group of isometries is compact semisimple) and symmetric spaces of noncompact type (if the group of isometries is noncompact semisimple). Symmetric spaces of compact and noncompact type are in some way dual to each other, and some of their properties can be carried from one type to the other.  An example of this is the study of totally geodesic submanifolds.  However, many properties are very different.  Noncompact symmetric spaces are diffeomorphic  to Euclidean spaces, and thus their topology is trivial, whereas in compact symmetric spaces topology does play a relevant role.  In fact, symmetric spaces of noncompact type are isometric to solvable Lie groups endowed with a left-invariant metric.  In our experience, this provides a wealth of examples of many interesting concepts, compared to their compact counterparts.

Our aim when studying homogeneous submanifolds is two-fold.  Firstly, we are interested in the classification (maybe under certain conditions) of homogeneous submanifolds of a given Riemannian manifold up to isometric congruence.  Usually we focus on the codimension one case, that is, homogeneous hypersurfaces, but we are also interested in higher codimension under some additional assumptions, for example when the group of isometries acts on the manifold polarly.  An isometric action is said to be polar if there is a submanifold that intersects all orbits orthogonally.  Such a submanifold is called a section of the polar action.  If the section is flat, then the action is called hyperpolar.  Polar actions take their name from polar coordinates, a concept that they generalize.  Sections are usually seen as sets of canonical forms~\cite{PT:tams}, as it is often the case that in symmetric spaces sections are precisely the Jordan canonical forms of matrix groups.

The second problem that we would like to address is the characterization of (certain classes of) homogeneous submanifolds.  It is obvious that homogeneity imposes restrictions on the geometry of a submanifold, and this in turn has implications on its shape.  The question is whether a particular property imposed on shape by homogeneity is specific of homogeneous submanifolds or, on the contrary, there might be other submanifolds having this property.  For example, homogeneous hypersurfaces have constant principal curvatures, and it is known in Euclidean and real hyperbolic spaces that this property characterizes homogeneous hypersurfaces.  However, this is not the case in spheres, as there are examples of hypersurfaces with constant principal curvatures that are not homogeneous.  We are particularly interested in isoparametric hypersurfaces, that is, hypersurfaces whose nearby equidistant hypersurfaces have constant mean curvature.  It is easy to see that homogeneous hypersurfaces are isoparametric, but we will see in this survey to what extent the converse is true.

Finally, we also study CPC submanifolds, that is, submanifolds whose principal curvatures, counted with their multiplicities, are independent of the unit normal vector.  This turns out to be an interesting notion related to several other properties whose study has recently attracted our attention.

This survey is organized as follows. In Section~\ref{sec:symmetric} we review the current definition, basic properties and types of Riemannian symmetric spaces, as well as the algebraic characterization of their totally geodesic submanifolds. Then, we deal more deeply with symmetric spaces of noncompact type in Section~\ref{sec:noncompact}, giving special relevance to the so-called Iwasawa decomposition of a noncompact semisimple Lie algebra.  This implies that a symmetric space of noncompact type is isometric to certain Lie group with a left-invariant Riemannian metric. The simplest symmetric spaces, apart from Euclidean spaces, are symmetric spaces of rank one, which include spaces of constant curvature. Rank one symmetric spaces of noncompact type are studied in Section~\ref{sec:rank1}, where we discuss different results regarding homogeneous hypersurfaces, isoparametric hypersurfaces and polar actions.  Finally, we study symmetric spaces of higher rank in Section~\ref{sec:higher_rank}.  A refinement of the Iwasawa decomposition is obtained in terms of parabolic subgroups in this section, and this is used to provide certain results in this setting, such as an extension method for submanifolds and isometric actions. Moreover, we report on what is known about polar actions in this context, and explain a method to study  homogeneous CPC submanifolds given by subgroups of the solvable part of the Iwasawa decomposition of the symmetric space.

\section{A quick review on symmetric spaces}\label{sec:symmetric}
In this section we include a short introduction to symmetric spaces. We first present the notion and first properties (\S\ref{subsec:notion}), then the different types of symmetric spaces (\S\ref{subsec:types}), and we conclude with an algebraic characterization of totally geodesic submanifolds (\S\ref{subsec:totally}). 

There are several references that the reader may like to consult to obtain further information on this topic. Probably, the most well-known and complete references are Helgason's book~\cite{Helgason} and Loos' books~\cite{Loos1, Loos2}. Eschenburg's survey~\cite{Eschenburg} and Ziller's notes~\cite{Ziller} are great references, especially for beginners. The books by Besse~\cite{Besse}, Kobayashi and Nomizu~\cite{Kobayashi}, O'Neill~\cite{ONeill} and Wolf~\cite{Wolf} also include nice chapters on symmetric spaces. In this section we mainly follow~\cite{Helgason} and~\cite{Ziller}.

\subsection{The notion and first properties of a symmetric space.}\label{subsec:notion}

In any connected Riemannian manifold $M$ we can consider normal neighborhoods around any point $p\in M$. If we take a geodesic ball $B_p(r)=\{q\in M: d(p,q)<r\}$, with $r$ small enough, as one of these neighborhoods, we can always consider a smooth map $\sigma_p\colon B_p(r)\to B_p(r)$ that sends each $q=\exp_p(v)$ to $\sigma_p(q)=\exp_p(-v)$, for $v\in T_pM$, $|v|<r$; hereafter, $\exp$ denotes the Riemannian exponential map. The map $\sigma_p$ is an involution, i.e.\ $\sigma_p^2=\id$, which is called geodesic reflection. If, for any $p\in M$, one can define $\sigma_p$ in the same way globally in $M$ and $\sigma_p$ is an isometry of $M$, then we say that $M$ is a \emph{(Riemannian) symmetric space}.

It follows easily from the definition that symmetric spaces are complete (since geodesics can be extended by using geodesic reflections) and homogeneous, that is, for any $p_1$, $p_2\in M$ there is an isometry $\varphi$ of $M$ mapping $p_1$ to $p_2$ (take $\varphi=\sigma_q$, where $q$ is the midpoint of a geodesic joining $p_1$ and $p_2$). A Riemannian manifold $M$ is homogeneous if and only if the group $\Isom(M)$ of isometries of $M$ acts transitively on $M$. Then, $M$ is diffeomorphic to a coset space $G/K$ endowed with certain differentiable structure. Here, $G$ can be taken as the connected component of the identity element of the isometry group of $M$, i.e.\ $G=\Isom^0(M)$, which still acts transitively on $M$ since $M$ is assumed to be connected, whereas $K=\{g\in G: g(o)=o\}$ is the isotropy group of some (arbitrary but fixed) base point $o\in M$. As the isometry group of any Riemannian manifold is a Lie group, then $G$ is also a Lie group, and $K$ turns out to be a compact Lie subgroup of $G$.

Let us define the involutive Lie group automorphism $s\colon G\to G$, $g\mapsto \sigma_o g \sigma_o$, which satisfies $G_s^0\subset K\subset G_s$, where $G_s=\{g\in G:s(g)=g\}$ and $G_s^0$ is the connected component of the identity. The differential $\theta=s_*\colon \g{g}\to\g{g}$ of $s$ is a Lie algebra automorphism called the \emph{Cartan involution} of the symmetric space (at the Lie algebra level). The isotropy Lie algebra $\g{k}$ is the eigenspace of $\theta$ with eigenvalue $1$. Let $\g{p}$ be the $(-1)$-eigenspace of~$\theta$. The eigenspace decomposition of $\theta$ then reads $\g{g}=\g{k}\oplus\g{p}$, which is called the \emph{Cartan decomposition}. 
Moreover, it easily follows that $[\g{k},\g{k}]\subset \g{k}$, $[\g{k},\g{p}]\subset\g{p}$ and $[\g{p},\g{p}]\subset\g{k}$. This implies, by the definition of the Killing form $\cal{B}$ of $\g{g}$ (recall: $\cal{B}(X,Y)=\tr(\ad(X)\circ\ad(Y))$ for $X$, $Y\in\g{g}$), that $\g{k}$ and $\g{p}$ are orthogonal subspaces with respect to $\cal{B}$. 

By considering the map $\phi\colon G\to M$, $g\mapsto g(o)$, one easily gets that its differential $\phi_{*e}$ at the identity induces a vector space isomorphism $\g{p}\cong T_o M$. The linearization of the isotropy action of $K$ on $M$, which turns out to be the orthogonal representation $K\to\mathsf{GL}(T_oM)$, $k\mapsto k_{*o}$, is then equivalent to the adjoint representation $K\to \mathsf{GL}(\g{p})$, $k\to \Ad(k)$. Each one of these is called the \emph{isotropy representation} of the symmetric space. 

\subsection{Types of symmetric spaces.}\label{subsec:types}
If the restriction of the isotropy representation of $M\cong G/K$ to the connected component of the identity of $K$ is irreducible, we say that the symmetric space $M$ is \emph{irreducible}. This turns out to be equivalent to the property that the universal cover $\widetilde{M}$ of $M$ (which is always a symmetric space) cannot be written as a nontrivial product of symmetric spaces, unless $\widetilde{M}$ is some Euclidean space~$\R^n$.

A symmetric space $M\cong G/K$ is said to be of \emph{compact type}, of \emph{noncompact type}, or of \emph{Euclidean type} if $\cal{B}\rvert_{\g{p}\times\g{p}}$, the restriction to $\g{p}$ of the Killing form $\cal{B}$ of $\g{g}$, is negative definite, positive definite, or identically zero, respectively. If $M$ is irreducible, then Schur's lemma implies that $\cal{B}\rvert_{\g{p}\times\g{p}}$ is a scalar multiple of the induced metric on $\g{p}\cong T_oM$ and, according to the sign of such scalar, $M$ falls into exactly one of the three possible types. It turns out that if $M$ is of compact type, then $G$ is a compact semisimple Lie group, and $M$ is compact and of nonnegative sectional curvature; if $M$ is of noncompact type, then $G$ is a noncompact real semisimple Lie group, and $M$ is noncompact (indeed, diffeomorphic to a Euclidean space) and with nonpositive sectional curvature; and if $M$ is of Euclidean type, its Riemannian universal cover is a Euclidean space~$\R^n$. Moreover, in general, the universal cover of a symmetric space $M$ splits as a product $\widetilde{M}=M_0\times M_+\times M_-$, where $M_0=\R^n$ is of Euclidean type, $M_+$ is of compact type, and $M_-$ is of noncompact type.

Symmetric spaces of compact and noncompact type are related via the notion of duality. Being more specific, there is a one-to-one correspondence between simply connected symmetric spaces of compact type and (necessarily simply connected) symmetric spaces of noncompact type. Moreover, dual symmetric spaces have equivalent isotropy representations and, therefore, irreducibility is preserved by duality. Without entering into details, the trick at the Lie algebra level to obtain the dual symmetric space is to change $\g{g}=\g{k}\oplus\g{p}$ by the new Lie algebra $\g{g}^*=\g{k}\oplus i\g{p}$, where $i=\sqrt{-1}$. In spite of the simplicity of this procedure, dual symmetric spaces have, of course, very different geometric and even topological properties. Examples of dual symmetric spaces are the following: 
\begin{enumerate}
	\item The round sphere $\mathbb{S}^n=\SO_{n+1}/\SO_n$ and the real hyperbolic space $\R \Hy^n=\SO^0_{1,n}/\SO_n$. 
	\item As a extension of the previous example, the projective spaces over the division algebras (other than $\R$) and their dual hyperbolic spaces: the complex spaces $\C \P^n=\SU_{n+1}/\Sg(\U_1\U_n)$ and $\C \Hy^n=\SU_{1,n}/\Sg(\U_1\U_n)$, the quaternionic spaces $\H \P^n=\Sp_{n+1}/\Sp_1\Sp_n$ and $\H \Hy^n=\Sp_{1,n}/\Sp_1\Sp_n$, and the Cayley planes $\mathbb{O} \P^2=\mathsf{F}_4/\Spin_9$ and $\mathbb{O} \Hy^2=\mathsf{F}_4^{-20}/\Spin_9$. The spaces in this and the previous item, jointly with the real projective spaces $\R \P^n$, constitute the so-called \emph{rank one symmetric spaces}.
	\item The oriented compact Grassmannian $\mathsf{G}^+_{p}(\R^{p+q})=\SO_{p+q}/\SO_p\SO_q$ of all oriented $p$-dimensional subspaces of $\R^{p+q}$, and the dual noncompact Grassmannian $\mathsf{G}_{p}(\R^{p,q})=\SO^0_{p,q}/\SO_p\SO_q$, which para\-met\-rizes all $p$-dimensional timelike subspaces of the semi-Euclidean space $\R^{p,q}$ of dimension $p+q$ and signature $(p,q)$. This example can be extended to complex and quaternionic Grassmannians.
	\item Any compact semisimple Lie group $G$ with bi-invariant metric, whose coset space is given by $(G\times G)/\Delta G$, and its noncompact dual symmetric space $G^\C/G$, where $G^\C$ is the complex semisimple Lie group given by the complexification of $G$. For instance, $\SU_n$ and $\SL_n(\C)/\SU_n$ are dual symmetric spaces.
	\item The space $\SU_n/\SO_n$ of all Lagrangian subspaces of $\R^{2n}$, and its noncompact dual space $\SL_n(\R)/\SO_n$ of all positive definite symmetric matrices of determinant~$1$.
\end{enumerate}
The whole list of simply connected, irreducible symmetric spaces can be found, for example, in~\cite[pp.~516, 518]{Helgason}.

\begin{remark}
In some cases above we have written $M=G/K$, where the action of $G$ on $M$ is not necessarily effective (i.e.\ not necessarily $G=\Isom^0(M)$).  However, in all cases such $G$-action is almost effective, that is, the ineffective kernel $\{g\in G:g(p)=p, \text{ for all } p\in M\}$ of the $G$-action on $M$ is a discrete subgroup of~$G$. Being more precise, one always considers a so-called symmetric pair $(G,K)$, where $K$ is compact, there is an involutive automorphism $s$ of $G$ such that $G_s^0\subset K\subset G_s$, and $G$ acts almost effectively on $M=G/K$. For example, the complex hyperbolic space $\C H^n$ is usually expressed as $\SU_{1,n}/\Sg(\U_1\U_n)$ instead of $(\mathsf{SU}_{1,n}/\mathbb{Z}_{n+1})/(\Sg(\U_1\U_n)/\mathbb{Z}_{n+1})$, in spite of the fact that $\SU_{1,n}$ has the cyclic group $\mathbb{Z}_{n+1}$ as ineffective kernel. This practice is common in the study of symmetric spaces for simplicity reasons, and because all Lie algebras involved remain the same. The symmetric pairs $(G,K)$ of compact type with $G=\Isom^0(M)$ can be found in~\cite[pp.~324-325]{Wang-Ziller}.
\end{remark}

\subsection{Totally geodesic submanifolds.}\label{subsec:totally}
Among different kinds of Riemannian submanifolds, the totally geodesic ones typically play an important role. This is particularly true~in the case of symmetric spaces. Indeed, although the classification of totally geodesic submanifolds in symmetric spaces is still an outstanding problem, these submanifolds are, intrinsically, also symmetric, and admit a neat algebraic characterization, which we recall~below.

A vector subspace $\g{s}$ of a Lie algebra $\g{g}$ is called a \emph{Lie triple system} if $[[X, Y],Z]\in\g{s}$ for any $X,Y,Z\in\g{s}$. Let now $\g{g}=\g{k}\oplus\g{p}$ be a Cartan decomposition of a symmetric space $M\cong G/K$, corresponding to a base point $o\in M$, as above. A fundamental result states that, if $\g{s}$ is a Lie triple system of $\g{g}$ contained in $\g{p}$, then $\exp_o(\g{s})$ is a totally geodesic submanifold of $M$, and it is intrinsically a symmetric space itself. And conversely, if $S$ is a totally geodesic submanifold of $M$, and $o\in S$, then $\g{s}:=T_oS\subset T_oM\cong\g{p}$ is a Lie triple system. In this situation, $\g{h}=[\g{s},\g{s}]\oplus\g{s}$ is the Cartan decomposition of the Lie algebra $\g{h}$ of the isometry group of the symmetric space $S$. Indeed, there is a one-to-one correspondence between $\theta$-invariant subalgebras of $\g{g}$ and Lie triple systems.

A consequence of the previous characterization is the fact that totally geodesic submanifolds of symmetric spaces are preserved under duality: if $\g{s}\subset\g{p}$ is a Lie triple system in $\g{g}=\g{k}\oplus\g{p}$, then $i\g{s}\subset i\g{p}$ is a Lie triple system in $\g{g}^*=\g{k}\oplus i\g{p}$. Moreover, if $\g{s}\subset \g{p}$ is a Lie triple system, then $S=\exp_o(\g{s})$ is an intrinsically flat submanifold if and only if $\g{s}$ is an abelian subspace of $\g{p}$ (i.e.\ $[\g{s},\g{s}]=0$). This follows  from the Gauss equation of submanifold geometry, the property that $S$ is totally geodesic, and the fact that the curvature tensor $R$ of a symmetric space at the base point $o$ is given by
\begin{equation}\label{eq:curvature}
R(X,Y)Z=-[[X,Y],Z], \qquad X,Y,Z\in T_oM\cong \g{p}.
\end{equation}
Thus, one defines the \emph{rank} of a symmetric space $M$ as the maximal dimension of a totally geodesic and flat submanifold of $M$ or, equivalently, the dimension of a maximal abelian subspace of $\g{p}$. Clearly, the rank is an invariant that is preserved under duality.

In spite of the above algebraic characterization of totally geodesic submanifolds of symmetric spaces and the fact that, by duality, one can restrict to symmetric spaces of compact type (or of noncompact type), the classification problem remains open. In particular, one does not know any efficient procedure to classify Lie triple systems in general. 

Totally geodesic submanifolds of rank one symmetric spaces are well known (see~\cite[~\S3]{Wolf:totally}). The case of rank two is much more involved, and has been addressed by Chen and Nagano~\cite{Chen-Nagano1, Chen-Nagano2} and Klein~\cite{Klein:tams, Klein:osaka}. Apart from these works, the subclass of the so-called reflective submanifolds has been completely classified by Leung~\cite{Leung:indiana, Leung:jdg}. A submanifold of a symmetric space $M$ is called \textit{reflective} if it is a connected component of the fixed point set of an involutive isometry of $M$; or, equivalently, if it is a totally geodesic submanifold such that the exponentiation of one (and hence all) normal space is also a totally geodesic submanifold. Finally, let us mention that the index of symmetric spaces (that is, the smallest possible codimension of a proper totally geodesic submanifold) has been recently investigated by Berndt and Olmos~\cite{BO:crelle, BO:blms, BO:jdg}, who proved, in particular, that the index of an irreducible symmetric space is bounded from below by the rank. 
Further information on totally geodesic submanifolds of symmetric spaces can be found in~\cite[\S11.1]{BCO}.

\section{Symmetric spaces of noncompact type and their Lie group model}\label{sec:noncompact}
In this section we focus on symmetric spaces of noncompact type. Our goal will be to explain that any symmetric space of noncompact type is isometric to a Lie group endowed with a left-invariant metric. The reader looking for more information or detailed proofs can consult, for instance, Eberlein's~\cite[Chapter~2]{Eberlein}, Helgason's~\cite[Chapter~VI]{Helgason} or Knapp's books~\cite[Chapter VI, \S4-5]{Knapp}. A nice survey that includes a detailed description of the space $\SL_n(\R)/\SO_n$ can be found in~\cite{Berndt:hyperpolar}. In this section we mainly follow \cite{Berndt:hyperpolar} and~\cite{Knapp}.

\begin{figure}
	\begin{tabular}{c@{\quad}c@{\quad}c}
		\includegraphics[width=0.3\textwidth]{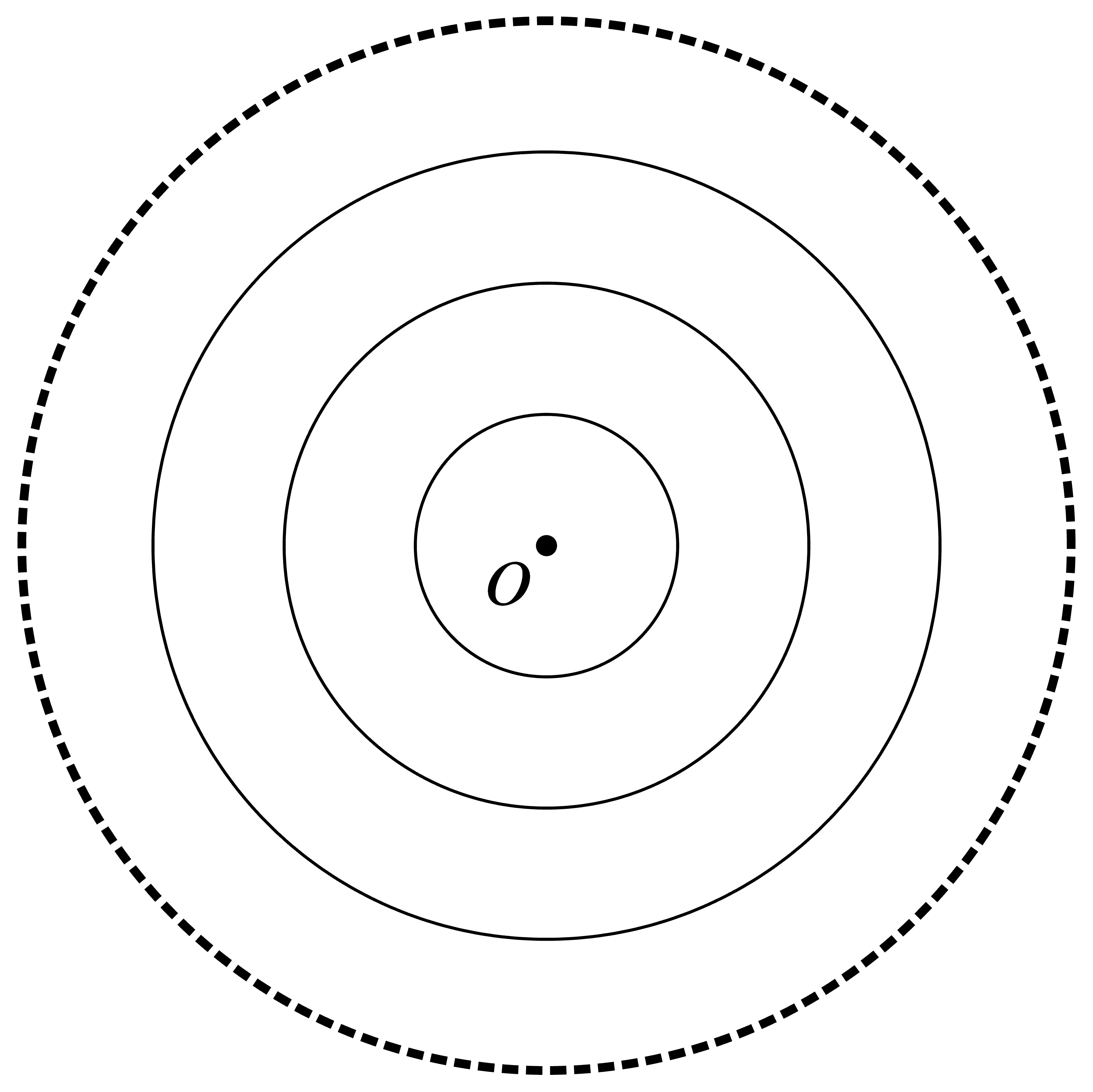}
	&
		\includegraphics[width=0.3\textwidth]{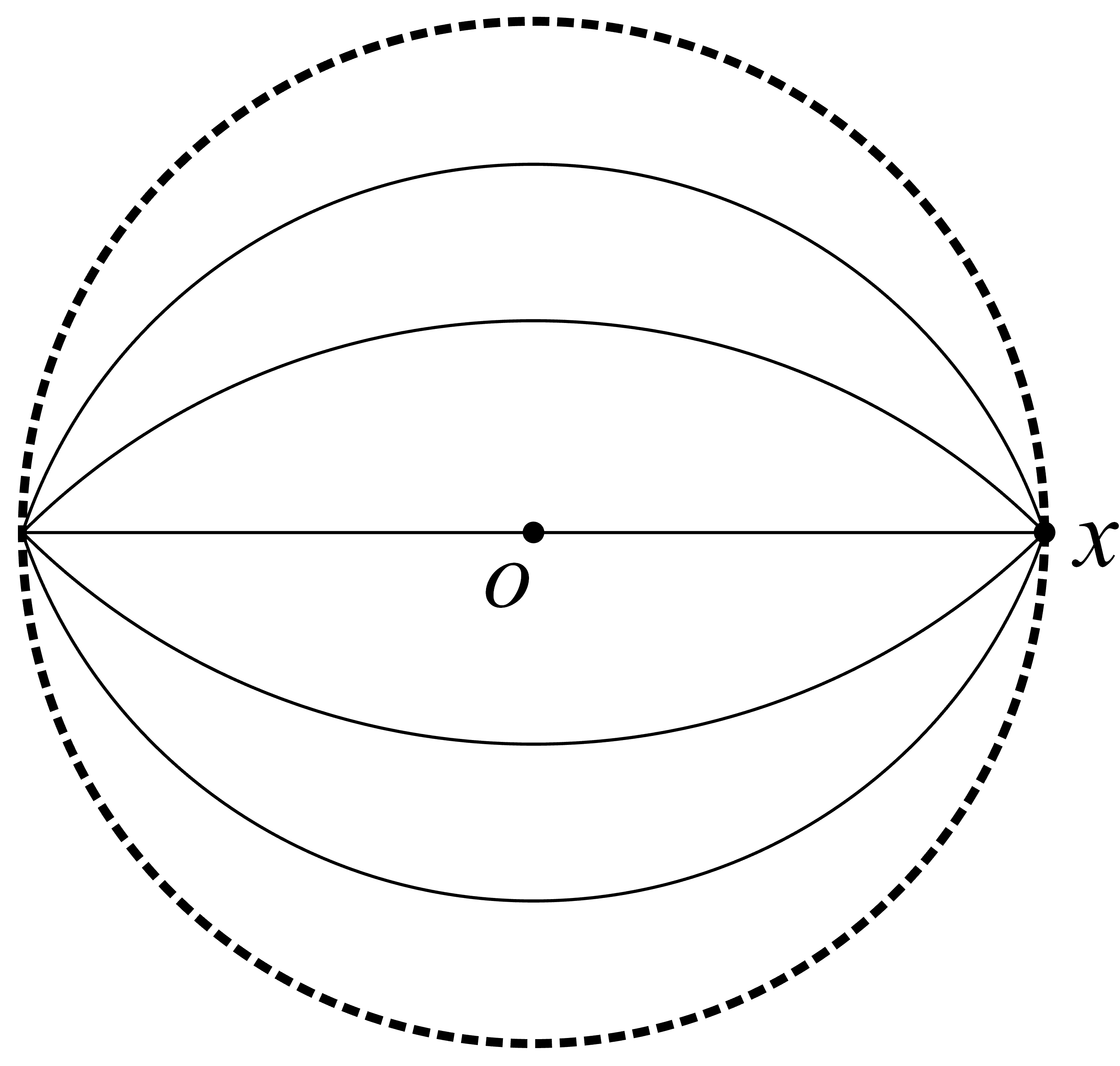}
	&
		\includegraphics[width=0.3\textwidth]{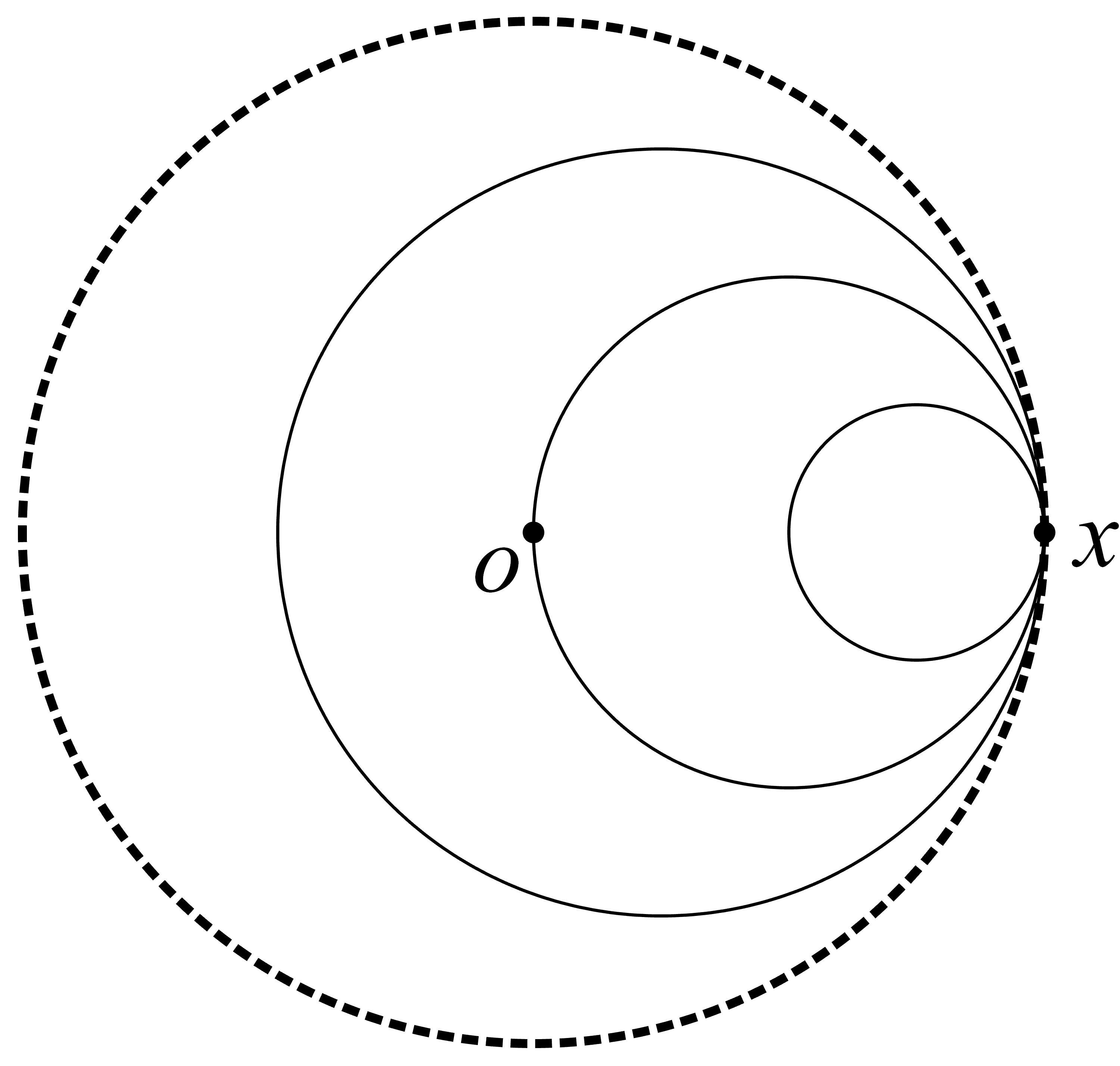}
	\end{tabular}
	\caption{Orbit foliations of the actions of the groups $K$, $A$ and $N$ on $\R \Hy^2$, respectively.}\label{fig:KAN}
\end{figure}

\begin{example}\label{ex:RH^2}
	The real hyperbolic plane $\R \Hy^2$ is the most basic example of symmetric space of noncompact type, and the only one of dimension at most two. It is well-known that (as any other symmetric space of noncompact type) it is diffeomorphic to an open ball, which gives rise, with an appropriate metric, to the Poincar\'e disk model for $\R \Hy^2$. Let us consider, however, the half-space model, by regarding $\R \Hy^2$ as the set $\{z\in \C: \Im z>0\}$ with metric $\langle\cdot,\cdot\rangle_{\R^2}/(\Im z)^2$. Then, the group $G=\SL_2(\R)$ acts transitively, almost effectively and by isometries on $\R \Hy^2$ via M\"obius transformations:
	\[
	\begin{pmatrix}
	a &b \\ c & d
	\end{pmatrix}\cdot z=\frac{az+d}{cz+d}.
	\]
	Then, the isotropy group $K$ at the base point $o=\sqrt{-1}$ is $\SO_2$, and hence $\R \Hy^2=\SL_2(\R)/\SO_2$. Moreover, any matrix in $\SL_2(\R)$ can be decomposed in a unique way~as
	\[
		\begin{pmatrix}
		a &b \\ c & d
		\end{pmatrix}
		=
		\begin{pmatrix}
		\cos s & \sin s \\ -\sin s & \cos s
		\end{pmatrix}
		\begin{pmatrix}
		\lambda &0 \\ 0 & \lambda^{-1}
		\end{pmatrix}
		\begin{pmatrix}
		1 &u \\ 0 & 1
		\end{pmatrix}, \qquad \text{where } s, u\in\R, \lambda>0.
	\]
	
	From an algebraic viewpoint, this decomposition turns out to encode some of the elements involved in the Gram-Schmidt process applied to the basis of $\R^2$ given by the column vectors of the left-hand side matrix: the orthogonal matrix is the transition matrix from the orthonormal basis produced by the method to the canonical basis, whereas the diagonal and upper triangular matrices contain the coefficients calculated in the process. 
	Moreover, the matrices on the right-hand side define three subgroups of $\SL_2(\R)$, namely $K=\SO_2$, the abelian subgroup $A$ of diagonal matrices, and the nilpotent subgroup $N$ of unipotent upper-triangular matrices. This so-called Iwasawa decomposition $G=KAN$ can be extended to any symmetric space of noncompact type, as we will soon explain.
	 
	From a geometric perspective, we can get insight into the groups involved in the decomposition by looking at their isometric actions on the hyperbolic plane (see~Figure~\ref{fig:KAN}). Thus, the $K$-action  fixes $o$ and the other orbits are geodesic spheres around $o$, the orbits of the $A$-action are a geodesic through $o$ and equidistant curves to such geodesic, while the $N$-action produces the horocycle foliation of $\R \Hy^2$ centered at one of the two points at infinity of the geodesic $A\cdot o$. This description of the actions make also intuitively clear the important fact that $\R \Hy^2\cong G/K$ is diffeomorphic to the subgroup $AN$ of $G$.
\end{example}

We now move on to the general setting. We start by describing some important decompositions of the Lie algebra of the isometry group (\S\ref{subsec:decompositions}), and then we present the Lie group model of a symmetric space of noncompact type~(\S\ref{subsec:model}).

\subsection{Root space and Iwasawa decompositions.}\label{subsec:decompositions}

Let $M\cong G/K$ be an arbitrary symmetric space of noncompact type. Then $\g{g}$ is a real semisimple Lie algebra, which implies that its Killing form $\cal{B}$ is nondegenerate. Indeed, the Cartan decomposition $\g{g}=\g{k}\oplus\g{p}$ is $\cal{B}$-orthogonal, $\cal{B}\rvert_{\g{k}\times\g{k}}$ is negative definite (due to the compactness of $K$), and $\cal{B}\rvert_{\g{p}\times\g{p}}$ is positive definite (since $M$ is of noncompact type). Hence, by reverting the sign on $\g{k}\times\g{k}$ or, equivalently, by defining $\cal{B}_\theta(X,Y)=-\cal{B}(\theta X, Y)$, for $X,Y\in\g{g}$, we have that $\cal{B}_\theta$ defines a positive definite inner product on $\g{g}$. It is easy to check that this inner product satisfies $\cal{B}_\theta(\ad (X)Y, Z)=-\cal{B}_\theta(Y, \ad(\theta X)Z)$, $X,Y,Z\in\g{g}$.

Let $\g{a}$ be a maximal abelian subspace of $\g{p}$. One can show that any two choices of $\g{a}$ are conjugate under the adjoint action of $K$ (similar to the fact that any two maximal abelian subalgebras of a compact Lie algebra are conjugate to each other). Moreover, by definition, the rank of $M\cong G/K$ is the dimension of $\g{a}$. 
For each $H\in\g{a}$, $X,Y\in\g{g}$, we have that
\[
\cal{B}_\theta(\ad(H)X,Y)=-\cal{B}_\theta(X,\ad(\theta H)Y)=\cal{B}_\theta(X,\ad(H)Y),
\]
which means that each operator $\ad(H)\in\End(\g{g})$ is self-adjoint with respect to $\cal{B}_\theta$. Moreover, if $H_1, H_2\in\g{a}$, then $[\ad(H_1),\ad(H_2)]=\ad[H_1,H_2]=0$, since $\ad\colon \g{g}\to \End(\g{g})$ is a Lie algebra homomorphism and $\g{a}$ is abelian. Thus, $\{\ad(H):H\in\g{a}\}$ constitutes a commuting family of self-adjoint endomorphisms of $\g{g}$. Therefore, they diagonalize simultaneously. Their common eigenspaces are called the \textit{restricted root spaces}, whereas their nonzero eigenvalues (which depend linearly on $H\in\g{a}$) are called the \textit{restricted roots} of $\g{g}$. In other words, if for each covector $\lambda\in\g{a}^*$ we define
\[
\g{g}_\lambda=\{X\in\g{g}: [H,X]=\lambda(H)X \text{ for all } H\in\g{a}\},
\]
then any $\g{g}_\lambda\neq 0$ is a restricted root space, and any $\lambda\neq 0$ such that $\g{g}_\lambda\neq 0$ is a restricted root. Note that $\g{g}_0$ is always nonzero, as $\g{a}\subset \g{g}_0$. If  $\Sigma=\{\lambda\in\g{a}^*:\lambda\neq 0,\,\g{g}_\lambda\neq 0\}$ denotes the set of restricted roots, then we have the following $\cal{B}_\theta$-orthogonal decomposition 
\begin{equation}\label{eq:root_space_decomposition}
\g{g}=\g{g}_0\oplus\biggl(\bigoplus_{\lambda\in\Sigma}\g{g}_\lambda\biggr),
\end{equation}
which is called the \emph{restricted root space decomposition} of $\g{g}$.

Observe that these definitions depend on the choice of $o\in M$ (or, equivalently, of a Cartan involution $\theta$ of $\g{g}$) and of the choice of the maximal abelian subspace $\g{a}$ of $\g{p}$. However, different choices of $o$ and $\g{a}$ give rise to decompositions that are conjugate under the adjoint action of $G$.  For simplicity, in this article we will not specify this dependence and we will also omit the adjective ``restricted". 

It is easy to check that the following properties are satisfied:
\begin{enumerate}[(i)]
	\item $[\g{g}_\lambda,\g{g}_\mu]\subset\g{g}_{\lambda+\mu}$, for any $\lambda$, $\mu\in\g{a}^*$.
	\item $\theta \g{g}_\lambda=\g{g}_{-\lambda}$ and, hence, $\lambda\in\Sigma$ if and only if $-\lambda\in \Sigma$.
	\item $\g{g}_0=\g{k}_0\oplus\g{a}$, where $\g{k}_0=\g{g}_0\cap\g{k}$ is the normalizer of $\g{a}$ in $\g{k}$.
\end{enumerate}
For each $\lambda\in \Sigma$, define $H_\lambda\in\g{a}$ by the relation $\cal{B}(H_\lambda,H)=\lambda(H)$, for all $H\in\g{a}$. Then we can introduce an inner product on $\g{a}^*$ by $\langle \lambda,\mu\rangle:=\cal{B}(H_\lambda,H_\mu)$. Thus, with a bit more work one can show that $\Sigma$ is an abstract root system in $\g{a}^*$, that is, it satisfies (cf.~\cite[\S{}II.5]{Knapp}):
\begin{enumerate}[(a)]
	\item $\g{a}^*=\spann\Sigma$,
	\item for $\alpha,\beta\in \Sigma$, the number
	$a_{\alpha\beta} =2\langle\alpha,\beta\rangle/
	\langle\alpha,\alpha\rangle$ is an integer,
	
	\item for $\alpha,\beta\in \Sigma$, we have
	$\beta-a_{\alpha\beta}\,\alpha\in \Sigma$.
\end{enumerate}
This system may be nonreduced, that is, there may exist $\lambda\in\Sigma$ such that $2\lambda\in \Sigma$.

Now we can define a positivity criterion on $\Sigma$ by declaring those roots that lie at one of the two half-spaces determined by a hyperplane in $\g{a}^*$ not containing any root to be positive. If $\Sigma^+$ denotes the set of positive roots, then $\Sigma=\Sigma^+\cup(-\Sigma^+)$. As is usual in the theory of root systems, one can consider a subset $\Pi\subset \Sigma^+$ of simple roots, that is, a basis of $\g{a}^*$ made of positive roots such that any $\lambda\in \Sigma$ is a linear combination of the roots in $\Pi$ where all coefficients are either nonnegative integers or nonpositive integers. Of course, the cardinality of $\Pi$ agrees with the dimension of $\g{a}$, i.e.\ with the rank of $G/K$. The set $\Pi$ of simple roots allows to construct the Dynkin diagram attached to the root system $\Sigma$, which is a graph whose nodes are the simple roots, and any two of them are joined by a simple (respectively, double, triple) edge whenever the angle between the corresponding roots is $2\pi/3$ (respectively, $3\pi/4$, $5\pi/6$); moreover, if the system is nonreduced, two collinear positive roots are drawn as two concentric nodes. 

Due to the properties of the root space decomposition, the subspace
\[
\g{n}=\bigoplus_{\lambda\in \Sigma^+}\g{g}_\lambda
\]
of $\g{g}$ is a nilpotent subalgebra of $\g{g}$. Moreover, $\g{a}\oplus\g{n}$ is a solvable subalgebra of $\g{g}$ such that $[\g{a}\oplus\g{n},\g{a}\oplus\g{n}]=\g{n}$. Any two choices of positivity criteria on $\Sigma$ give rise to isomorphic Dynkin diagrams and to nilpotent subalgebras $\g{n}$ that are conjugate by an element of $N_K(\g{a})$.

A fundamental result in what follows is the \emph{Iwasawa decomposition theorem}. At the Lie algebra level, it states that
\[
\g{g}=\g{k}\oplus\g{a}\oplus\g{n}
\]
is a vector space direct sum (but neither orthogonal direct sum nor semidirect product). Let us denote by $A$ and $N$ the connected Lie subgroups of $G$ with Lie algebras $\g{a}$ and $\g{n}$, respectively. Since $\g{a}$ normalizes $\g{n}$, the semidirect product $AN$ is the connected Lie subgroup of $G$ with Lie algebra $\g{a}\oplus\g{n}$. Then the Iwasawa decomposition theorem at the Lie group level states that the multiplication map
\[
K\times A\times N\to G,\qquad (k,a,n)\mapsto kan
\] 
is an analytic diffeomorphism, and the Lie groups $A$ and $N$ are simply connected. Indeed, as $A$ is abelian and $N$ is nilpotent, they are both diffeomorphic to Euclidean spaces~\cite[Theorem~1.127]{Knapp}. Hence, the semidirect product $AN$ is also diffeomorphic to a Euclidean~space.

\subsection{The solvable Lie group model.}\label{subsec:model}
Recall the smooth map $\phi\colon G\to M$, $h\mapsto h(o)$, from the end of \S\ref{subsec:notion} . The restriction $\phi\rvert_{AN}\colon AN\to M$ is injective; indeed, if $\phi(h)=\phi(h')$ with $h,h'\in AN$, then $h^{-1}h'(o)=o$, and hence $h^{-1}h'\in K\cap AN$, which, by the Iwasawa decomposition, implies that $h^{-1}h'=e$. It is also onto: if $p\in M$, then by the transitivity of $G$ there exists $h\in G$ such that $h(p)=o$, but using the Iwasawa decomposition we can write $h=kan$, with $k\in K$, $a\in A$, $n\in N$, and then $p=h^{-1}(o)=n^{-1}a^{-1}k^{-1}(o)=(an)^{-1}(o)$. Finally, $\phi\rvert_{AN}$ is a local diffeomorphism: as $\ker\phi_{*e}=\g{k}$, we have that $(\phi\rvert_{AN})_{*e}\colon \g{a}\oplus\g{n}\to T_oM$ is an isomorphism, and, by homogeneity, any other differential $(\phi\rvert_{AN})_{*h}$ is also bijective.

Therefore, $\phi\rvert_{AN}\colon AN\to M$ is a diffeomorphism. If we denote by $g$ the Riemannian metric on $M$, we can pull it back to obtain a Riemannian metric $(\phi\rvert_{AN})^*g$ on $AN$. Hence, we trivially have that $(M,g)$ and $(AN,(\phi\rvert_{AN})^*g)$ are isometric Riemannian manifolds.

Let now $h,h'\in AN\subset G$, and denote by $L_h$ the left multiplication by $h$ in $G$. Then 
\[
(h^{-1}\circ \phi\rvert_{AN}\circ L_h)(h')=h^{-1}(hh'(o))=h'(o)=\phi\rvert_{AN}(h'),
\]
from where we get $h^{-1}\circ \phi\rvert_{AN}\circ L_h=\phi\rvert_{AN}$ as maps from $AN$ to $M$. Hence, since $h^{-1}$ is an isometry of $(M,g)$, and using the previous equality, we have
\[
L_h^*(\phi\rvert_{AN})^*g=L_h^*(\phi\rvert_{AN})^*(h^{-1})^*g=(h^{-1}\circ \phi\rvert_{AN}\circ L_h)^*g=(\phi\rvert_{AN})^*g.
\]
This shows that $(\phi\rvert_{AN})^*g$ is a left-invariant metric on the Lie group $AN$. 

Altogether, we have seen that \textit{any symmetric space $M\cong G/K$ of noncompact type is isometric to a solvable Lie group $AN$ endowed with a left-invariant metric}. In particular, any symmetric space of noncompact type is diffeomorphic to a Euclidean space and, since it is nonpositively curved, it is a Hadamard manifold. This allows us to regard any of these spaces as an open Euclidean ball endowed with certain metric, as happens with the ball model of the real hyperbolic space. 

Moreover, it is sometimes useful to view a symmetric space of noncompact type $M$ as a dense, open subset of a bigger compact topological space $M\cup M(\infty)$ which, in this case, would be homeomorphic to a closed Euclidean ball. In order to do so, one defines an equivalence relation on the family of complete, unit-speed geodesics in $M$: if $\gamma$ and $\sigma$ are two of them, we declare them equivalent if they are asymptotic, that is, if $d(\gamma(t),\sigma(t))\leq C$, for certain constant $C$ and for all $t\geq 0$. Each equivalence class of asymptotic geodesics is called a \emph{point at infinity}, and the set $M(\infty)$ of all of them is the \emph{ideal boundary} of $M$. By endowing $M\cup M(\infty)$ with the so-called cone topology, $M\cup M(\infty)$ becomes homeomorphic to a closed Euclidean ball whose interior corresponds to $M$ and its boundary to $M(\infty)$. Two geodesics are asymptotic precisely when they converge to the same point in $M(\infty)$. We refer to~\cite[\S1.7]{Eberlein} for more details.

The Lie group model turns out to be a powerful tool for the study of submanifolds of symmetric spaces of noncompact type. The reason is that one can \textit{consider interesting types of submanifolds by looking at subgroups of $AN$} or, equivalently, at subalgebras of $\g{a}\oplus\g{n}$. For this reason, a good understanding of the root space decomposition is crucial. Of course, not every submanifold (even extrinsically homogeneous submanifold) of $M$ can be regarded as a Lie subgroup of $AN$, but very important types of examples arise in this way, sometimes combined with some additional constructions, as we will comment on in the following sections. In any case, if one wants to study submanifolds of $AN$ with particular geometric properties, one needs to have manageable expressions for the left-invariant metric on $AN$ and its Levi-Civita connection. We obtain the appropriate formulas below.

Let us denote by $\langle\cdot,\cdot\rangle_{AN}$ the inner product on $\g{a}\oplus\g{n}$ given by the left-invariant metric $(\phi\rvert_{AN})^*g$ on $AN$. Assume for the moment that $M$ is irreducible. Then, as mentioned in~\S\ref{subsec:notion}, the inner product $\phi^*g_o$ on $T_oM$ induced by the metric $g$ on $M$ is a scalar multiple of modified Killing form $\cal{B}_\theta$, i.e.\ $\phi^*g_o=k\cal{B}_\theta$, for some $k>0$. Let us define the inner product $\langle\cdot,\cdot\rangle:=k \cal{B}_\theta$ on $\g{g}$, and find the relation between $\langle\cdot,\cdot\rangle_{AN}$ and $\langle\cdot,\cdot\rangle$. Thus, if $X,Y\in\g{a}\oplus\g{n}$, and denoting orthogonal projections (with respect to $\cal{B}_\theta$) with subscripts, we have
\begin{equation}\label{eq:relation_inner} 
\begin{aligned}
\langle X,Y\rangle_{AN} &{}=
(\phi\rvert_{AN})^*g_o(X_\mathfrak{k}+X_\mathfrak{p},Y_\mathfrak{k}+Y_\mathfrak{p})=g_o(\phi_*X_\mathfrak{p},\phi_* Y_\mathfrak{p})= k \cal{B}_\theta(X_\mathfrak{p},Y_\mathfrak{p})
\\ 
&{}= k \cal{B}_\theta \biggl(\frac{1-\theta}{2}X,\frac{1-\theta}{2}Y\biggr)
=\frac{k}{4}\,\cal{B}_\theta(2X_\mathfrak{a}+X_\mathfrak{n}-\theta X_\mathfrak{n},2Y_\mathfrak{a}+Y_\mathfrak{n}-\theta Y_\mathfrak{n})
\\ 
&{}=\frac{k}{4}\left(4\cal{B}_\theta(X_\mathfrak{a},Y_\mathfrak{a})+\cal{B}_\theta(X_\mathfrak{n},Y_\mathfrak{n})+\cal{B}_\theta(\theta X_\mathfrak{n},\theta Y_\mathfrak{n})\right)\\
&{}= k\bigl(\cal{B}_\theta(X_\mathfrak{a},Y_\mathfrak{a})+\frac{1}{2}\cal{B}_\theta(X_\mathfrak{n},Y_\mathfrak{n})\bigr)
\\ 
&{}=
\langle X_\mathfrak{a},Y_\mathfrak{a}\rangle+\frac{1}{2}\langle X_\mathfrak{n},Y_\mathfrak{n}\rangle.
\end{aligned}
\end{equation}
If $M$ is reducible, one can adapt the argument (by defining $\langle\cdot,\cdot\rangle$ as a suitable multiple of $\cal{B}_\theta$ on each factor) to prove the same formula. Note that $\langle\cdot,\cdot\rangle$ inherits from $\cal{B}_\theta$ the property 
\begin{equation}
\label{eq:prop_B_theta}
\langle\ad (X)Y, Z\rangle=-\langle Y, \ad(\theta X)Z\rangle, \qquad \text{for } X,Y,Z\in\g{g}.
\end{equation}

Using Koszul formula, and relations~\eqref{eq:relation_inner} and~\eqref{eq:prop_B_theta}, one can obtain an important formula for the Levi-Civita connection $\nabla$ of the Lie group $AN$. Indeed, if $X,Y,Z\in\g{a}\oplus\g{n}$, and taking into account that $[\g{a}\oplus\g{n},\g{a}\oplus\g{n}]\subset\g{n}$, we have
\begin{equation}\label{eq:Levi-Civita}
\begin{aligned}
\langle \nabla_X Y,Z\rangle_{AN} 
&{}= \frac{1}{2}\bigl(\langle [X,Y],Z\rangle_{AN} - \langle [Y,Z],X\rangle_{AN} -\langle [X,Z],Y\rangle_{AN} \bigr)
\\
&{}= \frac{1}{4}\bigl( \langle [X,Y],Z\rangle - \langle [Y,Z],X\rangle -\langle [X,Z],Y\rangle \bigr)
\\ 
&{}=\frac{1}{4}\langle [X,Y]+[\theta X,Y]-[X,\theta Y],Z\rangle.
\end{aligned}
\end{equation}
Note that we started and finished with different inner products. Thus, in order to obtain an explicit formula for $\nabla_X Y$ one has to impose some restrictions on $X$ and $Y$. For example, if $X$ and $Y$ do not belong to the same root space, then $[\theta X,Y]$ and $[X,\theta Y]$ are orthogonal to $\g{a}$, whence in this case $2\nabla_X Y=\bigl([X,Y]+[\theta X,Y]-[X,\theta Y])_{\g{a}\oplus\g{n}}$.

\section{Submanifolds of rank one symmetric spaces}\label{sec:rank1}

In this section we review some results about certain important types of submanifolds in the rank one symmetric spaces of noncompact type. We start by describing these spaces in further detail (\S\ref{subsec:rank1}), and then we comment on homogeneous hypersurfaces (\S\ref{subsec:rank1_homogeneous}), isoparametric hypersurfaces (\S\ref{subsec:rank1_isoparametric}), and polar actions (\S\ref{subsec:rank1_polar}) on these spaces.

\subsection{Rank one symmetric spaces}\label{subsec:rank1}
As mentioned in \S\ref{subsec:types}, the symmetric spaces of noncompact type and rank one are the hyperbolic spaces $\mathbb{F} \Hy^n$, $n\geq 2$, over the distinct division algebras, $\mathbb{F}\in\{\R,\C,\H, \mathbb{O}\}$ ($n=2$ if $\mathbb{F}=\mathbb{O}$). We observe that $\C \Hy^1$, $\H \Hy^1$ and $\mathbb{O}\Hy^1$ are isometric (up to rescaling of the metric) to $\R \Hy^2$, $\R \Hy^4$ and $\R \Hy^8$, respectively. 

The isotropy representation (i.e.\ the adjoint action of $K$ on $\g{p}$) of rank one symmetric spaces is transitive on the unit sphere of $\g{p}$. Therefore, these Riemannian manifolds are not only homogeneous, but also \textit{isotropic}, which implies that they are \textit{two-point homogeneous}. Indeed, two-point homogeneous Riemannian manifolds are symmetric and, except for Euclidean spaces, have rank one~\cite{Szabo}, \cite[\S8.12]{Wolf}, and therefore the only noncompact examples (other than Euclidean spaces) are the symmetric spaces of noncompact type and rank one. Moreover, they are precisely the \textit{symmetric spaces of strictly negative sectional curvature} (even more, their sectional curvature is pinched between $c$ and $c/4$ for some $c<0$).

The root space decomposition~\eqref{eq:root_space_decomposition} of a symmetric space $M=\mathbb{F}\Hy^n$ of rank one is rather simple. One can show that $\Sigma=\{-\alpha,\alpha\}$ if $M=\R \Hy^n$, and $\Sigma=\{-2\alpha,-\alpha,\alpha,2\alpha\}$ otherwise. Thus the root space decomposition can be rewritten as
\[
\g{g}=\g{k}_0\oplus\g{a}\oplus\g{g}_{-2\alpha}\oplus\g{g}_{-\alpha}\oplus\g{g}_\alpha\oplus\g{g}_{2\alpha}.
\]
Of course, $\g{a}\cong\R$, $\g{g}_{-2\alpha}\cong\g{g}_{2\alpha}$ and $\g{g}_{-\alpha}\cong\g{g}_{\alpha}$. Note that the connected subgroup $K_0$ of $K$ with Lie algebra $\g{k}_0$ normalizes each one of the spaces in the decomposition. See Table~\ref{table:rank_one} for the explicit description of the group $K_0$ and the spaces in the decomposition. Determining all this information involves a few linear algebra computations; see~\cite[Chapter~2]{DV:dea} for the case of the complex hyperbolic space $\C \Hy^n$.

According to \S\ref{subsec:model}, a symmetric space $M\cong G/K$ of noncompact type is isometric to a Lie group $AN$ with a left-invariant metric. In the rank one setting, $A$ is one-dimensional, and, by declaring $\alpha$ as a positive root, $N$ can be taken to be the connected subgroup of $G$ with Lie algebra $\g{n}=\g{g}_\alpha\oplus\g{g}_{2\alpha}$.

\renewcommand\arraystretch{1.25}
\begin{table}
	\begin{tabular}[h!]{llllll}
		\hline
		Symmetric space & $G$ & $K$ & $K_0$ & $\g{g}_\alpha$ & $\g{g}_{2\alpha}$
		\\ \hline
		Real hyperbolic space $\R \Hy^n$ & $\SO^0_{1,n}$ & $\SO_n$ & $\SO_{n-1}$ & $\R^{n-1}$ & $0$
		\\ \hline 
		Complex hyperbolic space $\C \Hy^n$ & $\SU_{1,n}$ & $\Sg(\U_1\U_n)$ & $\U_{n-1}$ & $\C^{n-1}$ & $\R$
		\\ \hline 
		Quaternionic hyperbolic space $\H \Hy^n$ & $\Sp_{1,n}$ & $\Sp_1\Sp_n$ & $\Sp_1\Sp_{n-1}$ & $\H^{n-1}$ & $\R^3$
		\\ \hline 
		Cayley hyperbolic plane $\mathbb{O} \Hy^2$ & $\mathsf{F}_4^{-20}$ & $\Spin_9$ & $\Spin_{7}$ & $\mathbb{O}$ & $\R^7$
		\\  \hline
	\end{tabular}
	\bigskip
	\caption{Symmetric spaces of noncompact type and rank one}\label{table:rank_one}
\end{table}

The geometric interpretation of the groups involved in the Iwasawa decomposition of~$G$ is similar to that of $M=\R \Hy^2$, described in Example~\ref{ex:RH^2} and Figure~\ref{fig:KAN}. The action of the isotropy group $K$ on $M=\mathbb{F}\Hy^{n}$ has $o$ as a fixed point, and the other orbits are geodesic spheres around $o$. The action of $A$ gives rise to a geodesic through $o$ (since $\g{a}$ is a Lie triple system), and the other orbits are equidistant curves. Note that any geodesic curve, such as $A\cdot o$, determines two points at infinity; the choice of a positivity criterion on the set $\Sigma$ of roots (equivalently, choosing $\g{n}=\g{g}_\alpha\oplus\g{g}_{2\alpha}$ or $\g{n}=\g{g}_{-\alpha}\oplus\g{g}_{-2\alpha}$) is interpreted geometrically as selecting one of the two points at infinity determined by $A\cdot o$. Thus, the orbits of the $N$-action are \textit{horospheres} centered precisely at the point at infinity $x$ determined by the choice of $\g{n}$. In other words, if $\gamma$ is a unit-speed geodesic such that $A\cdot o=\{\gamma(t):t\in\R\}$ and converging to $x\in M(\infty)$, then the orbits of the $N$-action are the level sets of the Busemann function $f_\gamma(p):=\lim_{t\to\infty}(d(\gamma(t),p)-t)$. See~\cite[\S2.2]{DV:dea} and~\cite[\S1.10]{Eberlein} for details.

\begin{example}\label{ex:horospheres}
	We will illustrate the use of Formula~\eqref{eq:Levi-Civita} by calculating the extrinsic geometry of horospheres in $\mathbb{F}\Hy^n$. Via the Lie group model, the horosphere $N\cdot o$ is nothing but the Lie subgroup $N$ of $AN$. Thus, its tangent space at any $g\in N$ is given by the left-invariant fields of $\g{n}$ at $g$. If $B\in\g{a}$ satisfies $\langle B,B\rangle_{AN}=1$, then it defines a unit normal vector field on $N$. Hence, the shape operator $\cal{S}$ of $N$ with respect to $B$ is given by $\cal{S}X=-\nabla_X B$, for $X\in\g{n}$. If $X\in\g{g}_\lambda$, for $\lambda\in\{\alpha,2\alpha\}$, then by the comment following~\eqref{eq:Levi-Civita}, we have
	\[
	\cal{S}X=-\nabla_X B=-\frac{1}{2}\bigl([X,B]+[\theta X,B]-[X,\theta B]\bigr)_{\g{a}\oplus\g{n}}
	=
	\lambda(B)X, 
	\]
	where we have used the definition of the root space $\g{g}_\lambda$, the fact that $[\theta X,B]\in\g{g}_{-\lambda}$ is orthogonal to $\g{a}\oplus\g{n}$, and $\theta(B)=-B$.  Hence, $N\cdot o$ has two distinct constant principal curvatures, $\alpha(B)$ and $2\alpha(B)$, with respective principal curvature spaces  $\g{g}_\alpha$ and  $\g{g}_{2\alpha}$.
	Finally, note that all horospheres in $M$ are congruent to each other. Indeed, any two horosphere foliations are congruent by an element in $K$. Moreover, the geodesic $A\cdot o$ intersects all $N$-orbits and, since $A$ normalizes $N$, any two $N$-orbits are congruent under an element~of~$A$.	
\end{example}

The Lie group model of a rank one symmetric space contains some underlying additional structure that is often very helpful. Let us define a linear map $J\colon\g{g}_{2\alpha}\to\End(\g{g}_\alpha)$ by 
\[
\langle J_ZU,V\rangle_{AN}=\langle [U,V],Z\rangle_{AN}, \qquad \text{for all }U,V\in\g{g}_\alpha, \, Z\in\g{g}_{2\alpha},
\]
or, equivalently by~\eqref{eq:relation_inner} and~\eqref{eq:prop_B_theta}, $J_ZU:=[Z,\theta U]$. Then, up to rescaling of the metric of $M$ (and hence of $\langle\cdot,\cdot\rangle_{AN}$), the endomorphism $J_Z$ satisfies (see~\cite[Proposition~1.1]{Koranyi:adv})
\[
J_Z^2=-\langle Z,Z\rangle_{AN} \Id_{\g{g}_\alpha},\qquad \text{for all }Z\in\g{g}_{2\alpha}.
\]
Thus, the map $J$ induces a representation of the Clifford algebra $\Cl\bigl(\g{g}_{2\alpha},-\langle\cdot,\cdot\rangle_{AN}\bigr)$ on~$\g{g}_\alpha$ (see~\cite[Chapter~1]{LM} for more information on Clifford algebras and their representations). This converts $AN$ with the rescaled left-invariant metric into a so-called \emph{Damek-Ricci space}, and its nilpotent part $N$ into a \emph{generalized Heisenberg group}. These concepts were introduced by Damek and Ricci~\cite{Damek-Ricci:bams} and by Kaplan~\cite{Kaplan:tams}, respectively, and a comprehensive work for their study is~\cite{BTV}. Regarding rank one symmetric spaces of noncompact type as Damek-Ricci spaces has the advantage of allowing to use the power of the theory of Clifford modules to obtain more manageable formulas and more general arguments. For example, Formula~\eqref{eq:Levi-Civita} for the Levi-Civita connection of $AN$ adopts the form 
\begin{equation}\label{eq:Levi-Civita_DR}
{\nabla}_{aB+U+X}(bB+V+Y)
=\left(\frac{1}{2}\langle U,V\rangle_{AN}\!+\langle X,Y\rangle_{AN}\!\right)\!B-\frac{1}{2}\left(bU+J_XV+J_YU\right)
+\frac{1}{2}[U,V]-bX,
\end{equation}
where $a$, $b\in\R$, $U$, $V\in\g{g}_\alpha$, $X$, $Y\in\g{g}_{2\alpha}$.

In what follows we review some results about certain types of submanifolds and isometric actions on symmetric spaces of noncompact type and rank one.

\subsection{Homogeneous hypersurfaces}\label{subsec:rank1_homogeneous}
A submanifold $P$ of a Riemannian manifold $M$ is said to be \emph{(extrinsically) homogeneous} if for any $p, q\in P$ there exists an isometry $\varphi$ of the ambient manifold $M$ such that $\varphi(p)=q$ and $\varphi(P)=P$. Equivalently, $P$ is a homogeneous submanifold if it is an orbit of an isometric action on $M$, i.e.\ there exists a subgroup $H$ of $\Isom(M)$ such that $P=H\cdot p$ for some $p\in P$. Moreover, $P$ is embedded if and only if $H=\{\varphi\in \Isom(M):\varphi(P)=P\}$ is closed in $\Isom(M)$, which means that the associated isometric action is proper. From now on, \emph{isometric actions will be assumed to be proper}.

\begin{remark}
The collection of orbits of an isometric action is the standard example of a singular Riemannian foliation. A  \emph{singular Riemannian foliation} $\cal{F}$ on a Riemannian manifold $M$ is a decomposition of $M$ into connected, injectively immersed submanifolds $L\in\cal{F}$ (called leaves) such that they are locally equidistant to each other, and there is a collection of smooth vector fields on $M$ that spans all tangent spaces to all leaves; see~\cite{AB:book},~\cite{ABT:dga} for more information on this concept. Singular Riemannian foliations can have leaves of different dimensions: the ones of highest dimension are called regular, and the others are singular. Orbit foliations, that is, singular Riemannian foliations induced by isometric actions, are sometimes called homogeneous foliations.
\end{remark}

When an isometric action has codimension one orbits, then it is called a \emph{cohomogeneity one action}, and its codimension one orbits are \emph{homogeneous hypersurfaces}. The homogeneity property for hypersurfaces is a rather strong condition. This motivates the problem of classifying homogeneous hypersurfaces or, equivalently, cohomogeneity one actions up to orbit equivalence, in specific Riemannian manifolds, mainly in those with large isometry group. Such classification is known, for example, for Euclidean and real hyperbolic spaces (as a consequence of Segre's~\cite{Segre} and Cartan's~\cite{Cartan} works on isoparametric hypersurfaces, see~\S\ref{subsec:rank1_isoparametric}), irreducible symmetric spaces of compact type~\cite{Ko:tams}, and simply connected homogeneous $3$-manifolds with $4$-dimensional isometry group~\cite{DVM:arxiv}. Below we focus on the classification problem in symmetric spaces of noncompact type, and refer the reader to~\cite[\S6]{Berndt:notes} and~\cite[\S2.9.3 and Chapters 12-13]{BCO} for more information on cohomogeneity one actions.

As in any other Hadamard manifold, cohomogeneity one actions on symmetric spaces of noncompact type have at most one singular orbit~\cite[\S2]{BB:crelle} and no exceptional orbits~\cite[Corollary~1.3]{Ly:gd}. If there is one singular orbit, then the other orbits are homogeneous hypersurfaces which arise as distance tubes around the singular orbit. It there are no singular orbits, then all orbits are homogeneous hypersurfaces, and they define a regular Riemannian foliation of the ambient space. 

Cohomogeneity one actions on hyperbolic spaces have been investigated by Berndt, Br\"uck and Tamaru in a series of papers. Berndt and Br\"uck~\cite{BB:crelle} classified cohomogeneity one actions with a totally geodesic orbit on hyperbolic spaces $M=\mathbb{F} \Hy^n$:
\begin{theorem}\label{th:BB}
	Let $F$ be a totally geodesic singular orbit of a cohomogeneity one action on $\mathbb{F}\Hy^n$, $n\geq 2$. Then $F$ is congruent to one of the following totally geodesic submanifolds:
	\begin{itemize}
		\item in $\R \Hy^n$: $\{o\}$, $\R \Hy^1$, \dots, $\R \Hy^{n-1}$;
		\item in $\C \Hy^n$: $\{o\}$, $\C\Hy^1,\dots,\C \Hy^{n-1},\R\Hy^n$;
		\item in $\H \Hy^n$: $\{o\}$, $\H\Hy^1,\dots,\H \Hy^{n-1},\C\Hy^n$;	
		\item in $\mathbb{O} \Hy^2$: $\{o\}$, $\mathbb{O}\Hy^1,\H \Hy^2$.
	\end{itemize}
	Conversely, each of these totally geodesic submanifolds arises as the singular orbit of some cohomogeneity one action.
\end{theorem}
In particular, if $M\neq \R \Hy^n$, not every totally geodesic submanifold of $M$ defines homogeneous distance tubes. Moreover, it follows from Cartan's work~\cite{Cartan} that singular orbits of cohomogeneity one actions on $\R\Hy^n$ must be totally geodesic.

Berndt and Br\"uck~\cite{BB:crelle} also found examples of cohomogeneity one actions with a nontotally geodesic singular orbit (for $M\neq \R \Hy^n$). This important construction goes as follows.~Consider the Lie algebra $\g{a}\oplus\g{g}_\alpha\oplus\g{g}_{2\alpha}$ of $AN$. Take a subspace $\g{w}$  of $\g{g}_\alpha$ and define the subalgebra
\begin{equation}\label{eq:focal}
\g{s}_\g{w}:=\g{a}\oplus\g{w}\oplus\g{g}_{2\alpha}
\end{equation}
of $\g{a}\oplus\g{n}$. Assume that the orthogonal complement $\g{w}^\perp:=\g{g}_\alpha\ominus\g{w}$ of $\g{w}$ in $\g{g}_\alpha$ is such that
\begin{equation}\label{eq:condition_normalizer}
N^0_{K_0}(\g{w}) \text{ acts transitively on the unit sphere of } \g{w}^\perp
\end{equation}
where $N^0_{K_0}(\g{w})$ denotes the connected component of the identity of the normalizer of $\g{w}$ in $K_0$. Then, if $S_\g{w}$ is the connected subgroup of $AN$ with Lie algebra $\g{s}_\g{w}$, the group $N^0_{K_0}(\g{w}) S_\g{w}$ acts with cohomogeneity one on $M$, and with $S_\g{w}\cdot o$ as singular orbit if $\dim\g{w}^\perp\geq 2$. 

Berndt and Br\"uck proceeded to analyze which subspaces $\g{w}$ of $\g{g}_\alpha$ satisfy Condition~\eqref{eq:condition_normalizer}. In the case $M=\C\Hy^n$, they characterized this condition in terms of the so-called K\"ahler angles of $\g{w}^\perp$. Given any real subspace $V$ of a complex Euclidean space $(\R^{2k},J)$, where $J$ is a complex structure on $\R^{2k}$ (i.e.\ $J\in\g{so}(2k)$ and $J^2=-\Id$), the K\"ahler angle of a nonzero $v\in V$ with respect to $V$ is the angle $\varphi\in[0,\pi/2]$ between $Jv$ and $V$. When all unit $v\in V$ have the same K\"ahler angle $\varphi$ with respect to $V$, then we say that $V$ has \emph{constant K\"ahler angle} $\varphi$. For example, subspaces with constant K\"ahler angle $0$ or $\pi/2$ are precisely the complex and totally real subspaces, respectively. However, there are subspaces with any constant K\"ahler angle $\varphi\in(0,\pi/2)$; these can be classified, see~\cite[Proposition~7]{BB:crelle}. Now recall from Table~\ref{table:rank_one} that $\g{g}_\alpha\cong\C^{n-1}\cong(\R^{2n-2},J)$. Thus, it was proved in~\cite{BB:crelle} that $\g{w}\subset \g{g}_\alpha$ satisfies \eqref{eq:condition_normalizer} if and only if $\g{w}^\perp$ has constant K\"ahler angle $\varphi$ and $\dim\g{w}^\perp\geq 2$; moreover, the singular orbit $S_\g{w}\cdot o$ is nontotally geodesic whenever $\varphi\neq 0$.

In the case $M=\mathbb{O}\Hy^2$, by analyzing the $\Spin_7$-action on $\g{g}_\alpha\cong\R^8$, Berndt and Br\"uck proved that $\g{w}$ satisfies \eqref{eq:condition_normalizer} if and only if $\dim\g{w}\in\{0,1,2,4,5,6\}$, where only $\g{w}=0$ yields a totally geodesic singular orbit. Interestingly, the case $M=\H \Hy^n$ is much more involved and, indeed, it is still open. In~\cite{BB:crelle} it was proved that Condition~(7) in this case implies that $\g{w}^\perp$ has constant quaternionic K\"ahler angle (see~\S\ref{subsec:rank1_isoparametric}, after Theorem~\ref{th:isopar_DR}, for the definition), and several subspaces with this property were found (more examples were constructed in~\cite{DRDV:adv}). However, neither a classification of subspaces of $\H^{k}$, $k\geq 2$, with constant quaternionic K\"ahler angle, nor the equivalence between this property and Condition~\eqref{eq:condition_normalizer} are known. 

Regarding cohomogeneity one actions without singular orbits, Berndt and Tamaru~\cite{BT:jdg} proved (as a particular case of a more general result, cf.~\S\ref{subsec:cohom1}) that there are only two such actions on $M=\mathbb{F}\Hy^n$ up to orbit equivalence. One of these actions is that of the nilpotent part $N$ of the Iwasawa decomposition, giving rise to a \emph{horosphere foliation} (see Example~\ref{ex:horospheres}). The other one is given by the action of the connected subgroup $S$ of $AN$ with Lie algebra $\g{s}=\g{a}\oplus(\g{g}_\alpha\ominus\R U)\oplus\g{g}_{2\alpha}$, for any $U\in\g{g}_\alpha$; note that this corresponds to~\eqref{eq:focal} for the choice of a hyperplane $\g{w}$ in $\g{g}_{2\alpha}$. This $S$-action gives rise to the so-called \emph{solvable foliation} on a symmetric space of noncompact type and rank one.

Based on the results mentioned above, Berndt and Tamaru~\cite{BT:tams} were able to prove a structure result for cohomogeneity one actions on rank one symmetric spaces which states that each of these actions must be of one of the types described above.
\begin{theorem}\label{th:BT1}
	Let $M=\mathbb{F}\Hy^n$ be an symmetric space of noncompact type and rank one, and let $H$ act on $M$ with cohomogeneity one. Then one of the following statements holds:
	\begin{enumerate}[{\rm (1)}]
		\item The $H$-orbits form a regular Riemannian foliation on $M$ which is congruent to either a horosphere foliation or a solvable foliation.
		\item There exists exactly one singular $H$-orbit and one of the following two cases holds:
		\begin{enumerate}[{\rm (i)}]
			\item The singular $H$-orbit is one of the totally geodesic submanifolds in Theorem~\ref{th:BB}.
			\item The $H$-action is orbit equivalent to the action of $N^0_{K_0}(\g{w})S_\g{w}$, where $\g{w}$ is a subspace of $\g{g}_\alpha$ such that $N^0_{K_0}(\g{w})$ acts transitively on the unit sphere of $\g{w}^\perp$, and $S_\g{w}$ is the connected subgroup of $AN$ with Lie algebra $\g{s}_\g{w}=\g{a}\oplus\g{w}\oplus\g{g}_{2\alpha}$.
		\end{enumerate}
	\end{enumerate}
\end{theorem}
Combining this theorem with the results in~\cite{BB:crelle}, Berndt and Tamaru derived the classification of cohomogeneity one actions up to orbit equivalence on $\R \Hy^n$ and $\C\Hy^n$ for $n\geq 2$, and on $\H \Hy^2$ and $\mathbb{O}\Hy^2$. The classification on $\H \Hy^n$, $n\geq 3$, remains open.

\subsection{Isoparametric hypersurfaces} \label{subsec:rank1_isoparametric}
An immersed hypersurface $P$ in a Riemannian manifold $M$ is an \emph{isoparametric hypersurface} if, locally, $P$ and its nearby equidistant hypersurfaces have constant mean curvature. An \emph{isoparametric family of hypersurfaces} or \emph{isoparametric foliation (of codimension one)} is a singular Riemannian foliation such that its regular leaves are isoparametric hypersurfaces. These objects have been studied since the beginning of the 20th century and their investigation has therefore a long and interesting history. We refer to the excellent surveys~\cite{Chi:story} and~\cite{Th:survey} for a detailed account on this history.

Segre~\cite{Segre} classified isoparametric hypersurfaces in Euclidean spaces $\R^n$ by proving that they must be open subsets of affine hyperplanes $\R^{n-1}$, spheres $\mathbb{S}^{n-1}$ or generalized cylinders $\R^k\times \mathbb{S}^{n-k-1}$. Cartan~\cite{Cartan} proved that, in spaces of constant curvature, a hypersurface is isoparametric if and only if it has constant principal curvatures. Then, he classified such hypersurfaces in real hyperbolic spaces $\R \Hy^n$: the examples must be open subsets of totally geodesic $\R \Hy^{n-1}$ or their equidistant hypersurfaces, distance tubes around totally geodesic $\R \Hy^k$, $k\in\{0,\dots,n-2\}$, or horospheres. Thus, in spaces of nonpositive constant curvature, isoparametric hypersurfaces are open parts of homogeneous hypersurfaces. 

Observe that homogeneous hypersurfaces are isoparametric and have constant principal curvatures. However, none of the converse implications is true. In round spheres $\mathbb{S}^n$ there are inhomogeneous isoparametric hypersurfaces (with constant principal curvatures)~\cite{FKM}. In fact, the classification problem in spheres is much more involved; for more information, we refer the reader to some of the latest advances in the topic, such as~\cite{CCJ:annals,Chi:jdg,Imm:annals,Mi:annals,Siffert}. In spaces of nonconstant curvature, the problem becomes very complicated. Apart from the results we will review below, there is a classification on complex projective spaces $\C \P^n$, $n\neq 15$~\cite{DV:tams}, quaternionic projective spaces $\H \P^n$, $n\neq 7$~\cite{DVG:tohoku}, the product $\mathbb{S}^2\times\mathbb{S}^2$~\cite{Urbano}, and simply connected homogeneous $3$-manifolds with $4$-dimensional isometry group~\cite{DVM:arxiv}, such as the products $\mathbb{S}^2\times\R$, $\R \Hy^2\times\R$, the Heisenberg group $\textsf{Nil}_3$ or the Berger spheres. Interestingly, in all the cases mentioned so far (as well as in the rest of examples presented in this paper) an isoparametric hypersurface is always an open subset of a leaf of an isoparametric foliation of codimension one that fills the whole ambient space.

In spaces of nonconstant curvature, isoparametricity and constancy of the principal curvatures are two properties with no general theoretical relation. Berndt~\cite{Berndt:Hopf}, \cite{Berndt:HHn} classified curvature-adapted hypersurfaces with constant principal curvatures in complex and quaternionic hyperbolic spaces. Here,  \emph{curvature-adapted} means that the shape operator $\cal{S}$  and the normal Jacobi operator $R_\xi=R(\cdot,\xi)\xi$ of the hypersurface commute (hereafter $\xi$ is a unit normal smooth field on the hypersurface); hence both operators diagonalize simultaneously, which simplifies calculations involving the fundamental equations of submanifolds (Gauss, Codazzi...) and Jacobi fields adapted to the hypersurface (to calculate, for example, the extrinsic geometry of equidistant hypersurfaces or focal sets, cf.~\cite[\S10.2]{BCO}). In the complex case, a hypersurface in $\C\Hy^n$ is curvature-adapted if and only if it is \emph{Hopf}, that is, the Reeb vector field $J\xi$ is an eigenvector of the shape operator at every point, where $J$ is the K\"ahler structure of $\C\Hy^n$. It follows from Berndt's classifications that all curvature-adapted hypersurfaces with constant principal curvatures in $\C \Hy^n$ and $\H \Hy^n$ are open subsets of homogeneous hypersurfaces. However, not all homogeneous hypersurfaces described in~\S\ref{subsec:rank1_homogeneous} are curvature-adapted: only horospheres and homogeneous tubes around totally geodesic submanifolds have this property. Without the curvature-adaptedness condition, the study of hypersurfaces with constant principal curvatures is much more convoluted, and only some partial results for $\C \Hy^n$ are known; see~\cite{DRDV:indiana} for a recent advance, and~\cite{DV:dga} for a survey.

In view of the results mentioned above and the fact that a curvature-adapted hypersurface in a rank one symmetric space is isoparametric if and only if it has constant principal curvatures~\cite[Theorem~1.4]{GTY:japan}, it follows that a curvature-adapted isoparametric hypersurface in $\C \Hy^n$ or $\H \Hy^n$ is an open part of a homogeneous hypersurface. However, again, without the curvature-adaptedness condition, basically no other results regarding isoparametric hypersurfaces in our setting were known until a few years ago.

In~\cite{DRDV:mathz}, D\'iaz-Ramos and Dom\'inguez-V\'azquez constructed the first examples of inhomogeneous isoparametric hypersurfaces in a family of symmetric spaces of noncompact type, namely in complex hyperbolic spaces. Later, the authors generalized this result to Damek-Ricci spaces and, in particular, to the other symmetric spaces of noncompact type and rank one~\cite{DRDV:adv}. This construction, which we explain below, makes use of the basic idea of Berndt-Br\"uck cohomogeneity one actions described in~\S\ref{subsec:rank1_homogeneous}.

Given a symmetric space of noncompact type and rank one, $M=\mathbb{F}\Hy^n$, consider the subalgebra $\g{s}_\g{w}=\g{a}\oplus\g{w}\oplus\g{g}_{2\alpha}$ of $\g{a}\oplus\g{n}$ defined in~\eqref{eq:focal}, where now $\g{w}$ can be \emph{any} proper vector subspace of $\g{g}_\alpha$. Let $S_\g{w}$ be the connected subgroup of $AN$ with Lie algebra $\g{s}_\g{w}$. Using Formula~\eqref{eq:Levi-Civita_DR} it is not difficult to prove that $W_\g{w}:=S_\g{w}\cdot o$ is a minimal submanifold of $M$. Then, by introducing the notion of generalized K\"ahler angle (which we explain below) and using Jacobi field theory, D\'iaz-Ramos and Dom\'inguez-V\'azquez proved the following~\cite{DRDV:adv}:
\begin{theorem}\label{th:isopar_DR}
	The distance tubes around the minimal submanifold $W_\g{w}$ in a rank one symmetric space of noncompact type are isoparametric hypersurfaces, and have constant principal curvatures if and only if $\g{w}^\perp=\g{g}_\alpha\ominus\g{w}$ has constant generalized K\"ahler angle.
\end{theorem}
The concept of generalized K\"ahler angle extends both the K\"ahler angle and the quaternionic K\"ahler angle mentioned in~\S\ref{subsec:rank1_homogeneous}. Let $\g{z}\cong\R^m$ be a Euclidean space with inner product $\langle\cdot,\cdot\rangle$, and $\g{v}$  a Clifford module over $\Cl(\g{z},-\langle\cdot,\cdot\rangle)$. Consider $J\colon \g{z}\to \End(\g{v})$ the restriction to $\g{z}$ of the Clifford algebra representation. Recall from~\S\ref{subsec:rank1} that the rank one symmetric spaces of noncompact type have a naturally associated map $J$ as above, with $\g{v}=\g{g}_\alpha$ and $\g{z}=\g{g}_{2\alpha}$. Now let $V$ be a vector subspace of $\g{v}$. For each nonzero $v\in\g{v}$, consider the map
\[
F_v\colon \g{z}\to\R,\qquad Z\mapsto \langle (J_Zv)^V,(J_Zv)^V\rangle,
\] 
where $(\cdot)^V$ denotes orthogonal projection onto $V$. Observe that $F_v$ is a quadratic form on $\g{z}$, and its eigenvalues belong to the interval $[0,|v|^2]$. Hence, such eigenvalues are of the form $|v|^2\cos^2\varphi_i(v)$, $i=1,\dots, m=\dim\g{z}$, for certain angles $\varphi_i(v)\in[0,\pi/2]$. Then, one defines the \emph{generalized K\"ahler angle} of $v$ with respect to $V$ as the ordered $m$-tuple of angles $(\varphi_1(v),\dots,\varphi_m(v))$. We say that $V$ has \emph{constant generalized K\"ahler angle} if the $m$-tuple $(\varphi_1(v),\dots,\varphi_m(v))$ is independent of the nonzero $v\in V$. Note that, if $m =1$, we recover the notion of K\"ahler angle. The concept of quaternionic K\"ahler angle introduced in~\cite{BB:crelle} agrees with that of generalized K\"ahler angle in the case where $m=3$ and $\g{v}$ is a sum of equivalent irreducible $\Cl_3$-modules (i.e.\ $\g{v}$ is a quaternionic vector space).

Regarding complex or quaternionic hyperbolic spaces, $\C \Hy^n$ or $\H \Hy^n$, with $n\geq 3$, most real subspaces of $\g{g}_\alpha$ ($\cong \C^{n-1}$ or $\H^{n-1}$, respectively) have nonconstant generalized K\"ahler angle; e.g.\ the orthogonal sum of a complex and a totally real subspace in $\C^{n-1}$ does~not have constant K\"ahler angle. Thus, Theorem~\ref{th:isopar_DR} ensures the existence of  \emph{inhomogeneous isoparametric families of hypersurfaces with nonconstant principal curvatures in $\C \Hy^n$ and~$\H \Hy^n$,~$n\geq 3$}.

The case of the Cayley plane is even more interesting. As proved in~\cite{BB:crelle} and mentioned in~\S\ref{subsec:rank1_homogeneous}, if the subspace $\g{w}$ of $\g{g}_\alpha$ has dimension $3$, the tubes around $W_\g{w}$ are not homogeneous. However, any subspace of $\g{g}_\alpha\cong\mathbb{O}$ has constant generalized K\"ahler angle; in the case $\dim\g{w}=3$, the generalized K\"ahler angle of $\g{w}^\perp$ is $(0,0,0,0,\pi/2,\pi/2,\pi/2)$. Thus, the tubes around the corresponding $W_\g{w}$ constitute an \emph{inhomogeneous isoparametric family of hypersurfaces with constant principal curvatures in $\mathbb{O}\Hy^2$}. This is the only such example known in any symmetric space, apart from the FKM-examples in spheres~\cite{FKM}.

The homogeneous isoparametric foliations described in~\S\ref{subsec:rank1_homogeneous}, jointly with the inhomogeneous ones presented in this section, constitute an important family of examples which may encourage us to tackle the classification problem of isoparametric hypersurfaces in the rank one symmetric spaces of noncompact type. However, this is a much more complicated problem. Indeed, the only advance so far in this direction is the classification of isoparametric hypersurfaces in complex hyperbolic spaces obtained recently by the authors~\cite{DRDVSL:adv}. This constituted the first complete classification of isoparametric hypersurfaces in a complete family of symmetric spaces since Segre's~\cite{Segre} and Cartan's~\cite{Cartan} works in the~30s.
\begin{theorem}\label{th:isopCHn}
	Let $M$ be a connected real hypersurface in a complex hyperbolic space $\C\Hy^n$, $n\geq 2$. Then $M$ is isoparametric if and only if it is an open subset of one of the following:
	\begin{enumerate}[{\rm (i)}]
		\item A tube around a totally geodesic $\C\Hy^k$, $k\in\{0,\dots,n-1\}$.
		\item A tube around a totally geodesic $\R\Hy^n$.
		\item A horosphere.
		\item A leaf of a solvable foliation.
		\item A tube around a submanifold $W_\g{w}$, for some subspace $\g{w}$ of $\g{g}_\alpha$ with $\dim (\g{g}_\alpha\ominus\g{w})\geq 2$.
	\end{enumerate}
\end{theorem}
In particular, each isoparametric hypersurface in $\C\Hy^n$ is an open part of a complete, topologically closed leaf of a (globally defined) isoparametric foliation on $\C \Hy^n$. Either such foliation is regular (examples (iii) and (iv)) or has one singular orbit (examples (i), (ii) and (v)) which is minimal and homogeneous. Moreover, the homogeneous hypersurfaces in $\C \Hy^n$ are precisely those in examples (i) through (iv), and those in (v) with $\g{w}^\perp=\g{g}_\alpha\ominus\g{w}$ of constant K\"ahler angle. Thus, an isoparametric hypersurface in $\C \Hy^n$ is an open part of a homogeneous one if and only if it has constant principal curvatures.

The proof of Theorem~\ref{th:isopCHn} is rather involved. The starting point is to consider the Hopf map $\pi\colon \mathsf{AdS}^{2n+1}\to\C\Hy^n$ from the anti De Sitter spacetime $\mathsf{AdS}^{2n+1}$, and to prove that the preimage $\pi^{-1}(M)$ of a hypersurface $M$ in $\C \Hy^n$ is isoparametric (in a semi-Riemannian sense) if and only if $M$ is isoparametric. Since $\mathsf{AdS^{2n+1}}$ has constant curvature, $\pi^{-1}(M)$ is isoparametric precisely when it has constant principal curvatures. However, since $\pi^{-1}(M)$ is a Lorentzian hypersurface, its shape operator does not need to be diagonalizable. By analyzing each one of the four possible Jordan canonical forms for such shape operator, one can show (using elementary algebraic and geometric calculations) that three of them correspond to each one of the examples (i), (ii), (iii) above. Dealing with the fourth Jordan canonical form is much more convoluted, and requires delicate calculations with Jacobi fields and various geometric ideas. Finally, such Jordan form turns out to correspond with examples (iv) and (v) in Theorem~\ref{th:isopCHn}.

\subsection{Polar actions}\label{subsec:rank1_polar}
An isometric action on a Riemannian manifold $M$ is called \emph{polar} if there is a (a fortiori, totally geodesic) connected immersed submanifold $\Sigma$ of $M$ that intersects all orbits, and every such intersection is orthogonal. The submanifold $\Sigma$ is called a \emph{section} of the action; if $\Sigma$ is flat with respect to the induced metric, the action is called \emph{hyperpolar}. Cohomogeneity one actions constitute a particular case of hyperpolar actions.

The notion of polarity traces back at least to Dadok's classification~\cite{Dadok} of polar representations (equivalently, polar actions on round spheres): such polar actions coincide exactly with the isotropy representations of symmetric spaces, up to orbit equivalence. Later, polar actions have been studied mainly in the context of symmetric spaces of compact type: see~\cite{PT:jdg} (cf.~\cite{GK:agag}) for the classification in the rank one spaces, \cite{Ko:jdg} and~\cite{KL:blms} for the irreducible spaces of arbitrary rank, and~\cite{DR:polar} for a survey. For general manifolds, there are some topological and geometric structure results, see~\cite[Chapter~5]{AB:book} and~\cite{GZ:fixed}. Moreover, the notions of polar and hyperpolar action have been extended to the realm of singular Riemannian foliations by requiring the existence of sections through all points; see~\cite[Chapter~5]{AB:book}, \cite{ABT:dga}. Thus, homogeneous polar foliations are nothing but the orbit foliations of polar actions.

In symmetric spaces of noncompact type, very few results are known. The classification of polar actions on real hyperbolic spaces $\R\Hy^n$ follows from Wu's work~\cite{Wu:tams}. 

\begin{theorem}\label{th:polarRHn}
	A polar action on $\R\Hy^n$, $n\geq 2$, is orbit equivalent to one of the following:
	\begin{enumerate}[{\rm (i)}]
		\item  The action of $\SO_{1,k}\times Q$, where $k\in\{0,\dots,n-1\}$ and $Q$ is a compact subgroup of $\SO_{n-k}$ acting polarly on $\R^{n-k}$.
		\item The action of $N\times Q$, where $N$ is the nilpotent part of the Iwasawa decomposition of $\SO_{1,n}$, and $Q$ is a compact subgroup of $\SO_{n-k}$ acting polarly on $\R^{n-k}$.
	\end{enumerate}
\end{theorem}

The first classification result in a symmetric space of noncompact type and nonconstant curvature was achieved by Berndt and D\'iaz-Ramos~\cite{BDR:agag} for the complex hyperbolic plane $\C \Hy^2$. This classification consists of five cohomogeneity one actions and four cohomogeneity two actions, up to orbit equivalence. Interestingly, all of them can be characterized geometrically~\cite{DRDVVC:agag}. 
\begin{theorem}
	A submanifold of $\C \Hy^2$ is isoparametric if and only if it is an open part of a principal orbit of a polar action on $\C \Hy^2$.
\end{theorem}
Here, we refer to the notion of \emph{isoparametric submanifold} (of arbitrary codimension) given by Heintze, Liu and Olmos~\cite{HLO}, as a submanifold $P$ with flat normal bundle, whose parallel submanifolds have constant mean curvature in radial directions, and such that, for each $p\in P$, there is a totally geodesic submanifold $\Sigma_p$ such that $T_p\Sigma_p=\nu_p P$. Thus, an \emph{isoparametric foliation} (of arbitrary codimension) is a polar foliation whose regular leaves are isoparametric. The orbit foliations of polar actions constitute the main set of examples of isoparametric foliations.

Regarding cohomogeneity two polar actions on $\C\Hy^2$, one can additionally prove~\cite{DRDVVC:jga}:
\begin{theorem}
	A submanifold of $\C\Hy^2$ is an open subset of a principal orbit of a cohomogeneity two polar action if and only if it is a Lagrangian flat surface with parallel mean curvature. Moreover, such surfaces have parallel second fundamental form.
\end{theorem}

Coming back to the classification problem of polar actions on $\C \Hy^n$, the case $n=2$ was extended by D\'iaz-Ramos, Dom\'inguez-V\'azquez and Kollross to all dimensions~\cite{DRDVK:mathz}.

\begin{theorem}\label{th:polarCHn}
	A polar action on $\C \Hy^n$, $n\geq 2$, is orbit equivalent to the action of the connected subgroup $H$ of $\U_{1,n}$ with one of the following Lie algebras:
	\begin{enumerate}[{\rm (i)}]
\item $\g{h}=\g{q} \oplus \g{so}_{1,k}  \subset \g{u}_{n-k}
\oplus \g{su}_{1,k}$, $k\in\{0,\dots,n\}$, where the connected subgroup $Q$ of~$\U_{n-k}$ with Lie algebra $\g{q}$ acts polarly with a totally real section on~$\C^{n-k}$.

\item $\g{h}=\g{q} \oplus \g{b} \oplus \g{w} \oplus \g{g}_{2\alpha}\subset\g{su}_{1,n}$,
where $\g{b}$ is a subspace of~$\g{a}$, $\g{w}$ is a
subspace of~$\g{g}_{\alpha}$, and $\g{q}$ is a subalgebra of~$\g{k}_0$
which normalizes~$\g{w}$ and such that the connected subgroup of $\SU_{1,n}$
with Lie algebra $\g{q}$ acts polarly with a totally real section on
$\g{w}^\perp=\g{g}_\alpha\ominus\g{w}$.
	\end{enumerate}
\end{theorem}

In case (i), one $H$-orbit is a totally geodesic $\R \Hy^k$ and the other orbits are contained in the distance tubes around it. In item (ii), either $\g{b}=\g{a}$, in which case the orbit $H\cdot o$ contains the geodesic $A\cdot o$, or $\g{b}=0$, in which case $H\cdot o$ is contained in the horosphere $N\cdot o$. Moreover, in case (ii), any choice of real subspace $\g{w}\subset\g{g}_\alpha\cong\C^{n-1}$ gives rise to at least one polar action; the justification of this claim makes use of a decomposition theorem~\cite[\S2.3]{DRDVK:mathz} for real subspaces of a complex vector space as a orthogonal sum of subspaces of constant K\"ahler angle. Thus, whereas in $\C \Hy^2$ the moduli space of polar actions up to orbit equivalence is finite, in $\C \Hy^n$, $n\geq 3$, it is uncountable infinite.

\begin{remark}
	It is curious to observe that the orbit $H\cdot o$ corresponding to case (ii) in Theorem~\ref{th:polarCHn} with $\g{b}=\g{a}$ is precisely the singular leaf of the isoparametric foliations referred to in Theorems~\ref{th:isopar_DR} and~\ref{th:isopCHn}(v). In particular, it is a minimal submanifold, and the orbit foliation of the $H$-action constitutes a subfoliation of the isoparametric family of hypersurfaces given by the tubes around $H\cdot o$.
\end{remark}

\section{Submanifolds of symmetric spaces of arbitrary rank}\label{sec:higher_rank}
In this section we start by presenting some structure results for symmetric spaces of arbitrary rank, namely, their horospherical decomposition and the associated canonical extension procedure (\S\ref{subsec:horospherical}). Then we comment on the classification problem of cohomogeneity one and hyperpolar actions (\S\ref{subsec:cohom1}), and on a recent result on homogeneous CPC submanifolds~(\S\ref{subsec:CPC}).

\subsection{Horospherical decomposition and canonical extension.}\label{subsec:horospherical}
In this subsection we introduce these two important tools for the study of submanifolds in higher rank symmetric spaces. Further information can be found in~\cite[\S{}VII.7]{Knapp}, \cite[\S2.17]{Eberlein}, \cite[\S{}I.1]{BorelLi} and~\cite{DV:imrn}.

Let $M\cong G/K$ be a symmetric space of noncompact type. We follow the notation of Section~\ref{sec:noncompact}. Let $\Sigma$ be the set of roots of $M$, and $\Pi$ a set of simple roots, $|\Pi|=\rank M$. 

Let $\Phi$ be any subset of $\Pi$. Let $\Sigma_\Phi=\Sigma\cap\spann\Phi$ be the set of roots spanned by elements of $\Phi$, and $\Sigma_\Phi^+=\Sigma^+
\cap\spann\Phi$ the positive roots in $\Sigma_\Phi$. Then, we define
\[
\g{l}_\Phi=\g{g}_0\oplus\left(\bigoplus_{\lambda\in \Sigma_\Phi}\g{g}_\lambda \right), \qquad \g{n}_\Phi=\bigoplus_{\lambda\in \Sigma^+\setminus\Sigma^+_\Phi}\g{g}_\lambda,\qquad
\g{a}_\Phi=\bigcap_{\lambda\in\Phi}\ker \lambda,
\]
which are reductive, nilpotent and abelian subalgebras of $\g{g}$, respectively. Define also
\[
\g{m}_\Phi=\g{l}_\Phi\ominus\g{a}_\Phi,\qquad
\g{a}^\Phi=\g{a}\ominus\g{a}_\Phi=\bigoplus_{\lambda\in\Phi} \R H_\lambda.
\]
The subalgebra $\g{q}_\Phi=\g{l}_\Phi\oplus\g{n}_\Phi=\g{m}_\Phi\oplus\g{a}_\Phi\oplus\g{n}_\Phi$ is said to be the parabolic subalgebra of the real semisimple Lie algebra $\g{g}$ associated with the subset $\Phi\subset\Pi$. The decompositions $\g{q}_\Phi=\g{l}_\Phi\oplus\g{n}_\Phi$ and $\g{q}_\Phi=\g{m}_\Phi\oplus\g{a}_\Phi\oplus\g{n}_\Phi$ are known as the Chevalley and Langlands decompositions of $\g{q}_\Phi$, respectively.

\begin{remark} 
	By considering $L_\Phi$ as the centralizer of $\g{a}_\Phi$ in $G$, an $N_\Phi$ as the connected subgroup of $G$ with Lie algebra $\g{n}_\Phi$, one can define the parabolic subgroup $Q_\Phi=L_\Phi N_\Phi$ of $G$ associated with the subset $\Phi\subset\Pi$. Geometrically speaking, parabolic subgroups of $G$ are isotropy groups of points at infinity, i.e.\ $Q_\Phi=\{g\in G: g(x)=x\}$ for some $x\in M(\infty)$ (except for the case $\Phi=\Pi$, which gives rise to $Q_\Phi=G$). Thus, unlike points in $M$, isotropy groups of points at infinity are noncompact, and (except for $\rank M=1$) there are several (but finitely many, exactly $2^{\rank M}-1$) conjugacy classes of them.
\end{remark}

Consider the subspace 
\[
\g{b}_\Phi=\g{m}_\Phi\cap\g{p}=\g{a}^\Phi\oplus\biggl(\bigoplus_{\lambda\in\Sigma_\Phi^+}\g{p}_\lambda\biggr),
\]
where $\g{p}_\lambda=(1-\theta)\g{g}_\lambda$ is the orthogonal projection of $\g{g}_\lambda$ onto $\g{p}$. Then $\g{b}_\Phi$ is a Lie triple system (see~\S\ref{subsec:totally}) in $\g{p}$. We denote by $B_\Phi$ the corresponding totally geodesic submanifold of $M$ which, intrinsically, is a symmetric space of noncompact type and rank $|\Phi|$, and is known as the \emph{boundary component} of $M$ associated with the subset $\Phi\subset\Pi$. The Lie algebra of $\Isom(B_\Phi)$ is $\g{s}_\Phi:=[\g{b}_\Phi,\g{b}_\Phi]\oplus\g{b}_\Phi$. Thus, if $S_\Phi$ is the connected subgroup of $G$ with Lie algebra $\g{s}_\Phi$, then $B_\Phi=S_\Phi\cdot o$. It is not difficult to see that $B_\Phi$ can be regarded, under the isometry $M\cong AN$, as the connected subgroup of $AN$ with Lie algebra $\g{a}^\Phi\oplus\bigl(\bigoplus_{\lambda\in\Sigma_\Phi^+}\g{g}_\lambda\bigr)$.

The \emph{horospherical decomposition} theorem states that the map
\[
A_\Phi\times N_\Phi\times B_\Phi\to M,\qquad (a,n,p)\mapsto (an)(p),
\] 
is an analytic diffeomorphism, where $A_\Phi$ and $N_\Phi$ are the connected subgroups of $G$ with Lie algebras $\g{a}_\Phi$ and $\g{n}_\Phi$, respectively.

In other words, this result implies that the connected closed subgroup $A_\Phi N_\Phi$ of $AN$ acts isometrically and freely on $M$, and each $A_\Phi N_\Phi$-orbit intersects $B_\Phi$ exactly once. Moreover, such intersection is always orthogonal (see~\cite[Proposition~4.2]{BDRT:jdg}). Thus, \emph{the $A_\Phi N_\Phi$-action on $M$ is free and polar with section $B_\Phi$}. Moreover, as shown by Tamaru~\cite{Tamaru:math_ann}, \emph{all the orbits of the $A_\Phi N_\Phi$-action are Einstein solvmanifolds and minimal submanifolds of $M$} and, actually, they are mutually congruent by elements of $S_\Phi$. The $A_\Phi N_\Phi$-orbits are totally geodesic if and only if $\Phi$ and $\Pi\setminus\Phi$ are orthogonal sets of roots.

The reinterpretation of the horospherical decomposition as a free, polar action with minimal orbits gives rise to the so-called \emph{canonical extension method}, which was introduced in~\cite{BT:crelle} for cohomogeneity one actions, and generalized in~\cite{DV:imrn} to other types of actions, foliations and submanifolds. This method allows to extend such geometric objects from a boundary component $B_\Phi$ to the whole symmetric space $M$, that is, from symmetric spaces of lower rank to symmetric spaces of higher rank. And, more importantly, one can do so by preserving some important geometric properties. 

In order to formalize this, let $P$ be a submanifold of codimension $k$ in $B_\Phi$. Then
\[
A_\Phi N_\Phi\cdot P:=\{h(p):h\in A_\Phi N_\Phi,\, p\in P\}=\bigcup_{p\in P}A_\Phi N_\Phi\cdot p
\]
is a submanifold of codimension $k$ in $M$. The mean curvature vector field of $A_\Phi N_\Phi\cdot P$ is $A_\Phi N_\Phi$-equivariant and, along $P$, coincides with that of $P$. This implies that, if $P$ has parallel mean curvature, is minimal, has (globally) flat normal bundle, or is isoparametric, then $A_\Phi N_\Phi\cdot P$ has the same property. 

One can also extend singular Riemannian foliations from $B_\Phi$ to $M$ by extending their leaves as above. Thus, if $\cal{F}$ is a singular Riemannian foliation on $B_\Phi$ that is polar, hyperpolar, or isoparametric, then the extended foliation $A_\Phi N_\Phi\cdot \cal{F}=\{A_\Phi N_\Phi\cdot L:L\in\cal{F}\}$ has the same property. Moreover, if $\cal{F}$ is homogeneous, that is, if it is the orbit foliation of an isometric action of a subgroup $H\subset S_\Phi$ of isometries of $B_\Phi$, then the extended foliation $A_\Phi N_\Phi\cdot \cal{F}$ is the orbit foliation of the isometric action of $A_\Phi N_\Phi H$ on~$M$.

This technique plays an important role both in the construction of interesting types of submanifolds in symmetric spaces of higher rank, as well as in their classification. In~\cite{DV:imrn} it was used, for example, to extend the examples of inhomogeneous isoparametric hypersurfaces presented in~\S\ref{subsec:rank1_isoparametric} to symmetric spaces of higher rank and type $BC_r$, such as noncompact complex and quaternionic Grassmannians, or the complexified Cayley hyperbolic plane $\mathsf{E}_6^{-14}/\Spin_{10}\U_1$. Also, it was used to construct inhomogeneous isoparametric foliations of codimension higher than one on noncompact real Grassmannians, as well as new examples of polar but nonhyperpolar actions on spaces of rank higher than one.

\subsection{Cohomogeneity one, hyperpolar and polar actions.}\label{subsec:cohom1}
By the very definition of rank, cohomogeneity one and hyperpolar actions on rank one symmetric spaces constitute the same family of actions. In higher rank, there are hyperpolar actions of greater cohomogeneity. Moreover, in any rank, there are polar actions which are not hyperpolar; for example, $A_\Phi N_\Phi$ acts polarly but not hyperpolarly on $M$, whenever $\Phi\neq \emptyset$. However, the classification problem of any of these types of actions is widely open. 

The most general result regarding cohomogeneity one actions on symmetric spaces of noncompact type is due to Berndt and Tamaru~\cite{BT:crelle}:
\begin{theorem}\label{th:BT_higher}
	Let $M\cong G/K$ be an irreducible symmetric space of noncompact type, and let $H$ be a connected subgroup of $G$ acting on $M$ with cohomogeneity one. Then one of the following statements holds:
	\begin{enumerate}[{\rm (1)}]
		\item The orbits form a regular foliation on $M$ and the $H$-action is orbit equivalent to the action of the connected subgroup of $AN$ with one of the following Lie algebras:
		\begin{enumerate}[{\rm (i)}]
			\item  $(\g{a}\ominus\R X)\oplus\g{n}$ for some $X\in\g{a}$.
			\item $\g{a}\oplus(\g{n}\ominus\R U)$, where $U\in \g{g}_\lambda$, for some $\lambda\in \Pi$.
		\end{enumerate}
		\item There exists exactly one singular orbit and one of the following two cases holds:
		\begin{enumerate}[{\rm (i)}]
			\item $H$ is contained in a maximal proper reductive subgroup $L$ of $G$, the actions of $H$ and $L$ are orbit equivalent, and the singular orbit is totally geodesic.
			\item $H$ is contained in a maximal proper parabolic subgroup $Q_\Phi$ of $G$ and the $H$-action is orbit equivalent to one of the following:
			\begin{enumerate}[{\rm (a)}]
				\item The canonical extension of a cohomogeneity one action with a singular orbit on the boundary component $B_\Phi$ of $M$.
				\item The action of a group obtained by nilpotent construction.
			\end{enumerate}
		\end{enumerate}
	\end{enumerate}
\end{theorem}

Cohomogeneity one actions with no singular orbits, i.e.\ giving rise to homogeneous regular Riemannian foliations, were classified in~\cite{BT:jdg}; they correspond to case (1) in Theorem~\ref{th:BT_higher}. Note that they are induced by subgroups of $AN$. 

Cohomogeneity one actions with a totally geodesic singular orbit were classified in~\cite{BT:tohoku}, and they correspond to case (2)-(i) above. Interestingly, the associated orbit foliations arise as tubes around certain reflective submanifolds, except for a few exceptional cases. 

The case that is still open is that of cohomogeneity one actions with a nontotally geodesic singular orbit, case (2)-(ii) in Theorem~\ref{th:BT_higher}. The main difficulty has to do with the so-called \emph{nilpotent construction} method, which somehow extends the construction of cohomogeneity one actions with a nontotally geodesic singular orbit in rank one symmetric spaces (\S\ref{subsec:rank1_homogeneous}). We skip the explanation of the method here, and refer the reader to~\cite{BT:crelle} or to~\cite{BDV:tg}, where this method was investigated. In these papers one can also find the only complete classifications known so far on symmetric spaces of higher rank, namely on $\SL_3/\SO_3$, $\SL_3(\C)/\SU_3$, $\mathsf{G}_2^
2/\SO_4$, $\mathsf{G}_2^\C/\mathsf{G}_2$ and $\SO_{2,n}^
0/\SO_2\SO_n$, $n\geq 3$. 

In the more general setting of hyperpolar actions, the only known result is the classification of hyperpolar actions with no singular orbits on any symmetric space of noncompact type, up to orbit equivalence, due to Berndt, D\'iaz-Ramos and Tamaru~\cite{BDRT:jdg}. In other words, this result describes all hyperpolar homogeneous regular Riemannian foliations on symmetric spaces of noncompact type.

\begin{theorem}\label{th:BDRT}
	A hyperpolar action with no singular leaves on a symmetric space of noncompact type $M$ is orbit equivalent to the hyperpolar action of the connected subgroup of $AN$ with Lie algebra 
	\[
	(\g{a}\ominus V)\oplus\biggl(\g{n}\ominus\biggl(\bigoplus_{\lambda\in\Phi}\R X_\lambda\biggr)\biggr),
	\]
	where $\Phi\subset \Pi$ is any subset of mutually orthogonal  simple roots, and $V$ is any subspace~of~$\g{a}_\Phi$.
\end{theorem} 

Note that the condition $\langle\lambda,\mu\rangle=0$ for any $\lambda,\mu\in\Phi$ implies that the associated boundary component $B_\Phi$ is the Cartesian product of $|\Phi|$ symmetric spaces of rank one, $B_\Phi=\prod_{\lambda\in \Phi}\mathbb{F}_\lambda\Hy^{n_\lambda}$, $\mathbb{F}_\lambda\in\{\R,\C,\H,\mathbb{O}\}$. Thus, the intersection of the foliation described in Theorem~\ref{th:BDRT} with $B_\Phi$ is the product foliation of solvable foliations (cf.~\S\ref{subsec:rank1_homogeneous}) on each factor $\mathbb{F}_\lambda\Hy^{n_\lambda}$. The case $V=0$ corresponds to the canonical extension of such product foliation.

Regarding polar actions on symmetric spaces of rank higher than one, very little is known. Let us simply mention the classification of polar actions with a fixed point on any symmetric space by D\'iaz-Ramos and Kollross~\cite{DRK:dga}, and the investigation of polar actions by reductive subgroups due to Kollross~\cite{Ko:duality}.

\subsection{Homogeneous CPC submanifolds.}\label{subsec:CPC}
A submanifold $P$ of a Riemannian manifold $M$ will be called a \emph{CPC submanifold} if its principal curvatures, counted with multiplicities, are independent of the normal direction. In particular, a CPC submanifold is always austere (that is, the multiset of its principal curvatures is invariant under multiplication by~$-1$) and, hence, minimal. Although the terminology CPC comes from \emph{constant principal curvatures}, the property of being CPC is more restrictive than the one studied in~\cite{HOT91} (cf.~\cite[\S4.3]{BCO}). However, this notion is intimately related to cohomogeneity one actions. Indeed, if a cohomogeneity one action on a Riemannian manifold has a singular orbit, then the slice representation at any point of this orbit is transitive on the unit sphere of the normal space, which implies that all shape operators with respect to any unit normal vector are conjugate and, hence, the singular orbit is CPC.
The converse is not true in general. In fact, as mentioned after Theorem~\ref{th:BB}, there are totally geodesic (and, hence, CPC) submanifolds in the complex hyperbolic space whose distance tubes are not homogeneous. 

In real space forms, a submanifold is CPC if and only if the distance tubes around it are isoparametric hypersurfaces with constant principal curvatures. The necessity in this equivalence is no longer true in spaces of nonconstant curvature (a counterexample is the one mentioned in the previous paragraph, in view of Theorem~\ref{th:isopCHn}), but the sufficiency holds in any Riemannian manifold for submanifolds of codimension higher than one~\cite{GT:asian}.
Moreover, totally geodesic submanifolds are examples of CPC submanifolds. Thus, the study of CPC submanifolds encompasses important problems, such as the classifications of totally geodesic submanifolds, cohomogeneity one actions, and isoparametric hypersurfaces with constant principal curvatures. Let us also emphasize that the singular leaf of the inhomogeneous isoparametric family of hypersurfaces with constant principal curvatures on the Cayley hyperbolic plane described in~\S\ref{subsec:rank1_isoparametric} was, up to very recently, the only known example of a homogeneous, nontotally geodesic, CPC submanifold that is not an orbit of a cohomogeneity one action on a symmetric space of noncompact type.

In what follows we will report on the main results and ideas of a recent work by Berndt and Sanmart\'in-L\'opez~\cite{BS18} regarding CPC submanifods in irreducible symmetric spaces of noncompact type. One of the main goals of \cite{BS18} was precisely to provide a systematic approach to the construction of homogeneous, nontotally geodesic, CPC submanifolds, producing a large number of examples that are not orbits of cohomogeneity one actions. Another remarkable point is the introduction of an original and innovative technique based on the algebraic examination of the root system of symmetric spaces in order to calculate the shape operator of certain homogeneous submanifolds.

Let $M\cong G/K$ be an irreducible symmetric space of noncompact type; as usual, we follow the notation in~\S\ref{sec:noncompact}. Let $\alpha_0$, $\alpha_1\in \Pi$ be two simple roots connected by a single edge in the Dynkin diagram of the symmetric space $M$. Consider a Lie subalgebra $\g{s} = \g{a} \oplus (\g{n} \ominus V)$ of $\g{a}\oplus\g{n}$, where $V$ is a subspace of $\g{g}_{\alpha_0} \oplus \g{g}_{\alpha_1}$. This implies that $V = V_0 \oplus V_1$ with $V_k \subset \g{g}_{\alpha_k}$ for $k \in \{ 0, 1\}$. Let $S$ be the connected closed subgroup of $AN$ with Lie algebra $\g{s}$. In the following lines, we will explain the approach to the classification of the CPC submanifolds of the form $S \cdot o$. Moreover, in the final part of this section, we will see that with weaker hypotheses on $\g{s}$ we still achieve the same classification result.

The orbit $S \cdot o$ is a homogeneous submanifold and therefore it suffices to study its shape operator $\Ss$ at the point $o$. Since $\Ss_{\xi} X = - (\nabla_{X} \xi)^{\top}$, where $(\cdot)^{\top}$ denotes the orthogonal projection onto $\g{s}$, the idea is to analyze carefully the terms involved in the expression \eqref{eq:Levi-Civita} for the Levi-Civita connection of $M$. Let $\xi \in V$ be a unit normal vector to $S \cdot o$ and let $X \in \g{s}$ be a tangent vector to $S \cdot o$. First, assume that $X \in \g{a}$. Then
\[
[X, \xi] + [\theta X, \xi] - [X, \theta \xi] = - [X, \theta \xi] \in  \g{g}_{-\alpha_0} \oplus \g{g}_{-\alpha_1}.  
\]
Hence, $[X, \theta \xi]$ has trivial projection onto $\g{a} \oplus \g{n}$. Thus, $\Ss_\xi X=-(\nabla_{X} \xi)^\top = 0$ for any tangent vector $X \in \g{a}$ and any normal vector $\xi \in V$.
Now take $\xi \in V$ and $X \in \g{g}_{\lambda}^{\top}$ with $\lambda \in \Sigma^{+}$. Using~\eqref{eq:Levi-Civita} and some other considerations that we omit for the sake of simplicity, we obtain
\begin{equation}\label{equation:shape:operator}
\Ss_{\xi} X = - \dfrac{1}{2}  \left([X, \xi] -[X, \theta \xi]\right)^{\top}.
\end{equation}
Therefore, we deduce
\begin{equation}\label{equation:key:step}
\Ss_{\xi} X  \in (\g{g}_{\lambda + \alpha_0} \oplus \g{g}_{\lambda + \alpha_1}) \oplus (\g{g}_{\lambda - \alpha_0}^{\top} \oplus \g{g}_{\lambda - \alpha_1}^{\top}),
\end{equation}
for each $\xi \in V$ and each $X \in \g{g}_{\lambda}^{\top}$ with $\lambda \in \Sigma^{+}$. This shows that we need to understand how the shape operator $\Ss$ relates the different positive root spaces among them.

In order to clarify this situation, we introduce a generalization of the concept of $\alpha$-string \cite[p.~152]{Knapp}. For ${\alpha_0}, {\alpha_1} \in \Sigma$ and $\lambda \in \Sigma$ we define the $({\alpha_0}, {\alpha_1})$-string containing $\lambda$ as the set of elements in $\Sigma \cup \{0\}$ of the form $\lambda + n {\alpha_0} + m {\alpha_1}$ with $n,m \in \mathbb{Z}$. This allows to define an equivalence relation on $\Sigma^{+}$. We say that two roots $\lambda_1, \lambda_2 \in \Sigma^{+}$ are $(\alpha_0, \alpha_1)$-related if $\lambda_1 - \lambda_2 = n {\alpha_0} + m {\alpha_1}$ for some $n,m \in \mathbb{Z}$. Thus, the equivalence class $[\lambda]_{(\alpha_0, \alpha_1)}$ of the root $\lambda \in \Sigma^{+}$ consists of the elements which may be written as $\lambda +n {\alpha_0} + m {\alpha_1}$ for some $n,m \in \mathbb{Z}$. We will write  $[\lambda]$ for this equivalence class, taking into account that it depends on the roots $\alpha_0$ and $\alpha_1$ defining the string. Put $\Sigma^{+} / \sim$ for the set of equivalence classes. The family $\{[\lambda]\}_{\lambda \in \Delta^{+}}$ constitutes a partition of $\Sigma^{+}$. Using this notation, from~\eqref{equation:shape:operator}~and~\eqref{equation:key:step} we get that
\begin{equation}\label{eq:shape:class}
\Ss_{\xi} \left( \bigoplus_{\gamma \in [\lambda]} \g{g}_{\gamma}^{\top} \right) \subset \bigoplus_{\gamma \in [\lambda]} \g{g}_{\gamma}^{\top} \qquad\text{for all } \lambda \in \Sigma^{+}.
\end{equation}
This is the key point for studying if the orbit $S \cdot o$ is CPC. We will explain~\eqref{eq:shape:class} in words. For each $\lambda \in \Sigma^{+}$ the subspace $\bigoplus_{\gamma \in [\lambda]} \g{g}_{\gamma}^{\top}$ is a $\Ss_{\xi}$-invariant subspace of the tangent space $\g{s}$. Moreover, $S \cdot o$ is a CPC submanifold if and only if the eigenvalues of $\Ss_\xi$ are independent of the unit normal vector $\xi$ when restricted to each one of those invariant subspaces  $\bigoplus_{\gamma \in [\lambda]} \g{g}_{\gamma}^{\top}$, for every $\lambda \in \Sigma^{+}$. Thus it suffices to consider the orthogonal decomposition
\begin{equation}\label{invariant:decomposition}
\g{n} \ominus V  = \bigoplus_{[\lambda] \in \Sigma^{+} / \sim} \left(\bigoplus_{\gamma \in [\lambda]} \g{g}_{\gamma}^{\top} \right),
\end{equation}
and to study the shape operator when restricted to each one of these $\Ss_{\xi}$-invariant subspaces. Since $\alpha_0$ and $\alpha_1$ span an $A_2$ root system, then neither $2 \alpha_0$ nor $2 \alpha_1$ are roots. Hence, the $(\alpha_0, \alpha_1)$-string of $\alpha_0$ consists of the roots $\alpha_0$, $\alpha_1$ and $\alpha_0 + \alpha_1$. Thus, one of these subspaces is $\g{g}_{\alpha_0} \oplus \g{g}_{\alpha_1} \oplus \g{g}_{\alpha_0 + \alpha_1}$. This approach would be interesting if the rest of the subspaces respected some pattern and they could  be determined explicitly. The following result addresses both questions. Recall that ${\alpha_0}$ and ${\alpha_1}$ are simple roots connected by a single edge in the Dynkin diagram. We define the level of a positive root as the sum of the nonnegative coefficients of its expression with respect to the basis $\Pi$. Let $\lambda \in \Sigma^{+}$ be the root of minimum level in its $(\alpha_0, \alpha_1)$-string. Assume that it is not spanned by $\alpha_0$ and $\alpha_1$. Then, (taking indices modulo~$2$) we have:
\begin{enumerate}[{\rm (i)}]
	\item \label{structure:strings:i} If $\langle \lambda, \alpha_0 \rangle = 0 = \langle \lambda, \alpha_1 \rangle$, then $[\lambda] = \{\lambda\}$.
	\item If $|\alpha_k| \geq |\lambda|$ and $\langle \lambda, \alpha_k \rangle \neq 0$, then $[\lambda] =  \{\lambda, \lambda + {\alpha_k}, \lambda + {\alpha_{k}} + {\alpha_{k+1}}\}$. \label{structure:strings:ii}
	\item Otherwise, $[\lambda] =  \{\lambda, \lambda + {\alpha_k}, \lambda + {\alpha_k}  + {\alpha_{k+1}}, \lambda + 2{\alpha_k}, \lambda + 2 {\alpha_k} + {\alpha_{k+1}}, \lambda + 2 {\alpha_k} + 2{\alpha_{k+1}} \}$. \label{structure:strings:iii}
\end{enumerate}
The roots $\lambda$, $\alpha_0$ and $\alpha_1$ span a manageable subsystem and, roughly speaking, the proof follows from a case-by-case examination on the possible Dynkin diagrams for this subsystem. The CPC condition means that the eigenvalues of the shape operator do not depend on the normal vector when restricted to each one of the subspaces $\bigoplus_{\gamma \in [\lambda]} \g{g}_{\gamma}^{\top}$ in~\eqref{invariant:decomposition}, where $[\lambda]$ is one of the three possible types of strings above. 

If $\lambda$ is under the hypotheses of case~(\ref{structure:strings:i}), then $\g{g}_{\lambda}$ belongs to the $0$-eigenspace of the shape operator. This claim follows from~(\ref{equation:key:step}) and the fact that neither $\lambda + \alpha_k$ nor $\lambda - \alpha_k$ are roots for $k \in \{0,1\}$. 

We analyze case~(\ref{structure:strings:ii}) in order to give the key ideas for a nontrivial case. Let us start with some general considerations. For a fixed $l \in \{ 0,1 \}$, let $\gamma \in \Sigma^{+}$ be the root of minimum level in its $\alpha_l$-string, which consists of the roots $\gamma$ and $\gamma + \alpha_l$. Fix a normal unit vector $\xi_l \in V_{l}$ and define 
\begin{equation}\label{definition:phi}
\phi_{\xi_l} = |\alpha_l|^{-1} \ad(\xi_l) \qquad \text{and} \qquad \phi_{\theta \xi_l} = -|\alpha_l|^{-1} \ad(\theta\xi_l).
\end{equation}
These maps $\phi_{\xi_l}$ and $\phi_{\theta \xi_l}$ turn out to be inverse linear isometries in the sense that $\phi_{\theta\xi_l}\circ\phi_{\xi_l}\rvert_{\g{g}_\gamma}=\id_{\g{g}_\gamma}$ and $\phi_{\xi_l}\circ\phi_{\theta\xi_l}\rvert_{\g{g}_{\gamma+\alpha_l}}=\id_{\g{g}_{\gamma+\alpha_l}}$. Moreover, for each $X \in \g{g}_{\gamma}$ we have
\begin{equation}\label{easy:shape}
\nabla_X \xi_l = - \frac{|\alpha_l|}{2} {\phi_{\xi_l}(X)} \qquad \text{and} \qquad 
\nabla_{\phi_{\xi_l}(X)} \xi_l = - \frac{|\alpha_l|}{2} X.
\end{equation}
Let us come back to the study of case~(\ref{structure:strings:ii}). Write $\xi = \cos(\varphi) \xi_k + \sin(\varphi) \xi_{k+1}$ with $\varphi \in [0, \frac{\pi}{2}]$, $\xi_k \in V_k$ and $\xi_{k+1} \in V_{k+1}$. The following diagram may help to understand the situation. Note that the nodes represent root spaces and not roots.
\begin{equation*}
	\begin{tikzpicture}[scale=0.4]
	\draw(-15.0,0.) circle (0.5cm);
	\draw(-6.,0.) circle (0.5cm);
	\draw(3,0.) circle (0.5cm);
	\draw[-  triangle 45] (-14,0.3)-- (-7,0.3);
	\draw[-  triangle 45]  (-7,-0.3) -- (-14,-0.3) ;
	\draw[-  triangle 45] (-4,0.3)-- (2,0.3);
	\draw[-  triangle 45] (2,-0.3) -- (-4,-0.3);
	\begin{scriptsize}
	\draw[color=black] (-15, 1.1) node {$\g{g}_{\lambda}$};
	\draw[color=black] (-6, 1.1) node {$\g{g}_{\lambda + \alpha_k}$};
	\draw[color=black] (3, 1.1) node {$\g{g}_{\lambda + \alpha_k + \alpha_{k+1}}$};
	\draw[color=black] (-10.5, 0.8) node {$\phi_{\xi_k}$};
	\draw[color=black] (-10.5, -1.1) node {$\phi_{\theta \xi_k}$};
	\draw[color=black] (-1.5, 0.8) node {$\phi_{\xi_{k+1}}$};
	\draw[color=black] (-1.5, -1.1) node {$\phi_{\theta \xi_{k+1}}$};
	\end{scriptsize}
	\end{tikzpicture}
\end{equation*}
Take a unit tangent vector $X \in \g{g}_{\lambda}$. Using $\Ss_{\xi} X = - (\nabla_{X} \xi)^{\top}$ and~\eqref{easy:shape} for the pair $(\gamma, \alpha_l) \in \{(\lambda, \alpha_k),(\lambda + \alpha_k, \alpha_{k+1})\}$, we can see that the $3$-dimensional vector space spanned by $X, \phi_{\xi_k}(X),(\phi_{\xi_{k+1}} \circ \phi_{\xi_k})(X)$ is $\Ss_\xi$-invariant. The matrix representation of $\Ss_{\xi}$ restricted to $\g{g}_\lambda\oplus\g{g}_{\lambda+\alpha_k}\oplus\g{g}_{\lambda+\alpha_k+\alpha_{k+1}}$ is then given by $\dim (\g{g}_{\lambda})$ blocks of the form
\begin{equation}\label{cpc:matrix}
\frac{|{\alpha_0}|}{2} 
\left(\begin{array}{ccc} 0& \cos(\varphi) & 0 \\   \cos(\varphi) & 0 & \sin(\varphi)  \\  0& \sin(\varphi) & 0  
\end{array}\right),
\end{equation}
with respect to the decomposition $\g{g}_{ \lambda} \oplus \phi_{\xi_k}(\g{g}_{\lambda}) \oplus (\phi_{\xi_{k+1}} \circ \phi_{\xi_k})(\g{g}_{\lambda})$. The eigenvalues of the above matrix are $0$ and $\pm \frac{|\alpha_0|}{2}$. They do not depend on $\varphi$. It is also important to note that the nonzero principal curvatures depend on the length of the root $\alpha_0$. 

Case~(\ref{structure:strings:iii}) is slightly more difficult than the one we have just studied. Roughly speaking, it is necessary to generalize~(\ref{easy:shape}) having in mind the following diagram:
\begin{equation*}
	\begin{tikzpicture}[scale=0.4]
	\draw(-15.0,0.) circle (0.5cm);
	\draw(-6.,0.) circle (0.5cm);
	\draw(6,0.) circle (0.5cm);
	\draw(0,3.5) circle (0.5cm);
	\draw(0,-3.5) circle (0.5cm);
	\draw(15.,0.) circle (0.5cm);
	\draw[-  triangle 45] (7,0.) -- (14,0);
	\draw[-  triangle 45] (-14,0.0)-- (-7,0.0);
	\draw[-  triangle 45] (-5,0.35)-- (-1,3.15);
	\draw[triangle 45 -] (5,0.35) -- (1,3.15);
	\draw[-  triangle 45] (-5,-0.35)-- (-1,-3.15);
	\draw[triangle 45 -] (5,-0.35)-- (1,-3.15);
	\begin{scriptsize}
	\draw[color=black] (-10.5, 0.5) node {$\ad(\xi_k)$};
	\draw[color=black] (10.5, 0.5) node {$\ad(\xi_{k+1})$};
	\draw[color=black] (-3.75, 2.5) node {$\ad(\xi_k)$};
	\draw[color=black] (-4.2, -2.5) node {$\ad(\xi_{k+1})$};
	\draw[color=black] (4, 2.5) node {$\ad(\xi_{k+1})$};
	\draw[color=black] (3.75, -2.5) node {$\ad(\xi_k)$};
	\draw[color=black] (-15, 1.1) node {$\g{g}_{\lambda}$};
	\draw[color=black] (-6.2, 1.1) node {$\g{g}_{\lambda+\alpha_k}$};
	\draw[color=black] (6.5, 1.1) node {$\g{g}_{\lambda+2\alpha_k + \alpha_{k+1}}$};
	\draw[color=black] (0,4.3) node {$\g{g}_{\lambda+2\alpha_k}$};
	\draw[color=black] (0,-4.3) node {$\g{g}_{\lambda +  \alpha_k + \alpha_{k+1}}$};
	\draw[color=black] (15.3, 1.1) node {$\g{g}_{\lambda +2\alpha_k + 2 \alpha_{k+1}}$};
	\end{scriptsize}
	\end{tikzpicture}
\end{equation*}
The principal curvatures of the shape operator do not depend on the normal vector when restricted to subspaces induced by strings of type~(\ref{structure:strings:iii}). 

Thus, the problem can be reduced to studying the shape operator when restricted to 
\[
\g{g}_{\alpha_0}^\top \oplus \g{g}_{\alpha_1}^\top \oplus \g{g}_{\alpha_0 + \alpha_1}.
\]
In other words, one needs to study CPC submanifolds in a symmetric space with Dynkin diagram of type $A_2$, which would conclude the classification result. 

However, as mentioned above, we can state a more general result concerning CPC submanifolds. Denote by $\Pi^{\prime}$ the set of simple roots $\alpha \in \Pi$ with $2 \alpha \notin \Sigma$. Note that there is at most one simple root in $\Pi$ that does not belong to $\Pi^{\prime}$, and this happens when the restricted root system of $M$ is of type $BC_r$. Consider a Lie algebra $\g{s} = \g{a} \oplus (\g{n} \ominus V)$ with $V \subset \bigoplus_{\alpha \in \Pi'} \g{g}_{\alpha}$. This implies that $V = \bigoplus_{\alpha \in \Psi} V_{\alpha}$ for some set $\Psi \subset \Pi^{\prime}$. Similar ideas to those that led us to~(\ref{cpc:matrix}) allow to deduce that for each root $\alpha \in \Psi$ the nonzero eigenvalues of $\Ss_{\xi_{\alpha}}$ are proportional to the length of $\alpha$. Then,  if $\Psi$ contains roots $\alpha$ and $\beta$ of different lengths, it follows that the shape operators $\Ss_{\xi_{\alpha}}$ and $\Ss_{\xi_{\beta}}$ have different eigenvalues, which implies that $S \cdot o$ is not CPC. 

Moreover, assume that $\Psi$ contains at least three roots. Then $\Psi$ has at least two orthogonal roots, say $\alpha_0$ and $\alpha_1$. We will explain briefly why this cannot lead to a CPC submanifold $S \cdot o$. The key point is to find a positive root $\lambda \in \Sigma^{+}$ with nontrivial $\alpha_{k}$-string but trivial $\alpha_{k+1}$-string, for some $k \in \{0,1\}$ and indices modulo 2. According to~(\ref{easy:shape}), there must exist a tangent vector $X \in \g{g}_{\lambda} \oplus \g{g}_{\lambda + \alpha_k}$ such that $\Ss_{\xi_k} X = \mu X$, for a unit normal $\xi_k \in V_{\alpha_k}$ and some $\mu\neq 0$. However, from~(\ref{equation:key:step}) we deduce that $\Ss_{\xi_{k+1}} X = 0$ for a unit normal $\xi_{k+1} \in V_{\alpha_{k+1}}$. Thus, if we take a normal unit vector $\xi = \cos(\varphi) \xi_k + \sin(\varphi) \xi_{k+1}$ for $\varphi \in [0, 2 \pi]$, then we get $\Ss_{\xi} X = \cos(\varphi) \mu X$. Thus $S \cdot o$ cannot be CPC since we have an infinite family of different principal curvatures. 

Finally, if $V$ is contained in a single root space $\g{g}_{\alpha}$, $\alpha\in\Pi'$, then the $S$-action on $M$ is the canonical extension of a cohomogeneity one action with a totally geodesic orbit on the boundary component $B_{\{\alpha\}}\cong\mathbb{R}\Hy^n$ (see~\S\ref{subsec:horospherical}). Hence, if $\dim V\geq 2$,  $S\cdot o$ is a singular orbit of a cohomogeneity one action, and then CPC; if $\dim V=1$, $S\cdot o$ is the only minimal orbit of an action as in Theorem~\ref{th:BT_higher}(1-ii), which also happens to be CPC. Altogether, we can state the main result of~\cite{BS18}:
\begin{theorem}\label{th:BS}
	Let $\g{s} = \g{a} \oplus (\g{n} \ominus V)$ be a subalgebra of $\g{a} \oplus \g{n}$ with $V \subset \bigoplus_{\alpha \in \Pi'} \g{g}_{\alpha}$. Let $S$ be the connected closed subgroup of $AN$ with Lie algebra $\g{s}$. Then the orbit $S \cdot o$ is a CPC submanifold of $M \cong G/K$ if and only if one of the following statements holds:
	\begin{itemize}
		\item[{\rm (I)}] There exists a simple root $\lambda \in \Pi'$ with $V \subset \g{g}_{\lambda}$. \label{main:simple:examples}
		\item[{\rm (II)}] There exist two nonorthogonal simple roots $\alpha_0,\alpha_1 \in \Pi'$ with $|\alpha_0| = |\alpha_1|$ and subspaces $V_0 \subset \g{g}_{\alpha_0}$ and $V_1 \subset \g{g}_{\alpha_1}$ such that $V = V_0 \oplus V_1$ and one of the following conditions holds: \label{main:new:examples}
		\begin{itemize}
			\item[{\rm (i)}] $V_0 \oplus V_1 = \g{g}_{\alpha_0} \oplus \g{g}_{\alpha_1}$; \label{singular:orbit}
			\item[{\rm (ii)}] $V_0$ and $V_1$ are isomorphic to $\mathbb{R}$ and $V_0 \oplus V_1$ is a proper subset of $\g{g}_{\alpha_0} \oplus \g{g}_{\alpha_1}$; 
			\item[{\rm (iii)}] $V_0$ and $V_1$ are isomorphic to $\mathbb{C}$, $V_0 \oplus V_1$ is a proper subset of $\g{g}_{\alpha_0} \oplus \g{g}_{\alpha_1}$ and there exists $T \in \g{k}_0$ such that $\ad(T)$ defines complex structures on $V_0$ and $V_1$ and vanishes on $[V_0, V_1]$; \label{main:complex}
			\item[{\rm (iv)}] $V_0$ and $V_1$ are isomorphic to $\mathbb{H}$, $V_0 \oplus V_1$ is a proper subset of $\g{g}_{\alpha_0} \oplus \g{g}_{\alpha_1}$ and there exists a subset $\g{l} \subset \g{k}_0$ such that $\ad(\g{l})$ defines quaternionic structures on $V_0$ and $V_1$ and vanishes on $[V_0, V_1]$. \label{main:quaternionic}
		\end{itemize}
	\end{itemize}
	Moreover, only the submanifolds given by {\rm (I)} and {\rm (II)(i)} can appear as singular orbits of cohomogeneity one actions.
\end{theorem} 
\begin{remark}
	One may ask whether this result is still true if  $V$ is a subspace of the sum of root spaces corresponding to the roots in $\Pi$ (instead of $\Pi'$). However, this seems to be a more difficult problem. Indeed, it includes, in particular, the classification problem of CPC submanifolds of the type $S\cdot o$ in the quaternionic hyperbolic spaces $\H \Hy^n$, $n\geq 3$, which turns out to be equivalent to the classification of subspaces of $\g{g}_\alpha\cong \H^{n-1}$ with constant quaternionic K\"ahler angle. As we mentioned in~\S\ref{subsec:rank1_homogeneous}, this is nowadays an open problem.
\end{remark}

\section{Open problems}
We include a list of open problems related to the research presented above.
\begin{enumerate}
	\item In view of the exposition in~\S\ref{subsec:rank1_homogeneous}, classify homogeneous hypersurfaces in quaternionic hyperbolic spaces $\H \Hy^n$, $n\geq 3$. More generally, classify the real subspaces $V$ of $\H^{n-1}$ with constant quaternionic K\"ahler angle, and, for each one of them, determine if there is a subgroup of $\Sp_1\Sp_{n-1}$ that acts transitively on the unit sphere of $V$.
	\item Make further progress in the classification problem (mentioned in~\S\ref{subsec:rank1_isoparametric}) of hypersurfaces with constant principal curvatures in the complex hyperbolic spaces $\C \Hy^n$. This is a very difficult problem that lacks powerful ideas and techniques, apart from a clever combination of the information provided by the Codazzi and Gauss equations of submanifolds.
	\item Classify curvature-adapted hypersurfaces with constant principal curvatures in the Cayley hyperbolic plane $\mathbb{O}\Hy^2$.
	\item Initiate the investigation of non-curvature-adapted hypersurfaces with constant principal curvatures in $\H \Hy^n$ and $\mathbb{O} \Hy^2$.
	\item Is there any nonisoparametric hypersurface with constant principal curvatures in a symmetric space? In a general Riemannian setting, one can probably construct such a hypersurface for some specific ambient metric, but, even in this case, no concrete example is known to us.
	\item Prove that compact embedded hypersurfaces with constant mean curvature in rank one symmetric spaces of noncompact type must be geodesic spheres. Any idea towards the solution of this problem is likely to have profound implications and many applications to other problems (such as the symmetry of solutions to the so-called overdetermined boundary value problems, see e.g.~\cite{DVPSE}).
	\item Obtain a better understanding of the nilpotent construction method for cohomogeneity one actions, mentioned in Theorem~\ref{th:BT_higher}(2-ii-b). This seems to be a crucial step towards the solution of the classification problem of cohomogeneity one actions on irreducible symmetric spaces of noncompact type. Another approach may come from generalizing Theorem~\ref{th:BS} to the study of CPC submanifolds arising from arbitrary subgroups of the solvable part of the Iwasawa decomposition.
	\item Construct, if possible, new examples of inhomogeneous isoparametric hypersurfaces in symmetric spaces of noncompact type and rank higher than one. All the examples known so far arise as canonical extensions of isoparametric hypersurfaces in rank one symmetric spaces.
	\item Construct new examples of inhomogeneous isoparametric submanifolds of codimension higher than one on symmetric spaces of noncompact type. The only known examples are those appearing in Wu's classification for real hyperbolic spaces~\cite{Wu:tams}, and their canonical extensions to noncompact real Grassmannians. Such submanifolds are leaves of non-hyperpolar isoparametric foliations. Is there any example of an inhomogeneous hyperpolar foliation of codimension higher than one on a symmetric space of noncompact type, unlike the compact case~\cite{Christ} (cf.~\cite{Ly:gafa})?
	\item Make progress in the classification problem of totally geodesic submanifolds, mentioned in~\S\ref{subsec:totally}. This problem seems nowadays infeasible in full generality. However, with the algebraic method explained in \S\ref{subsec:CPC} we are able to calculate very efficiently the shape operator of many homogeneous submanifolds. These ideas may help to obtain some classification result in certain higher rank symmetric spaces.
\end{enumerate}

\end{document}